\documentclass{elsarticle}
\usepackage[utf8]{inputenc}
\usepackage{lineno,hyperref}
\modulolinenumbers[5]

\makeatletter
\def\ps@pprintTitle{%
  \let\@oddhead\@empty
  \let\@evenhead\@empty
  \let\@oddfoot\@empty
  \let\@evenfoot\@oddfoot
}
\makeatother

\usepackage{natbib}
\usepackage{graphicx}
\usepackage{multirow}
\usepackage{amsfonts}
\usepackage{hyperref}

\usepackage{amsmath}

\usepackage{caption}
\usepackage{subcaption}

\usepackage{xcolor}

\usepackage{comment}

\newcommand{\cut}[1]{}

\usepackage[ruled,vlined,linesnumbered]{algorithm2e}

\usepackage{hyperref}
\hypersetup{
     colorlinks = true,
     linkcolor = green!60!black,
     anchorcolor = green!50!black,
     citecolor = green!60!black,
     filecolor = green!60!black,
     urlcolor = green!60!black
     }
\usepackage[nameinlink,noabbrev]{cleveref}

\usepackage[normalem]{ulem}

\newcommand{\mathC}{\mathbb{C}}

\journal{Computer Physics Communications}

\newcommand{\rev}[1]{{\color{black}#1}}

\begin{document}

\begin{frontmatter}

\title{Coarsest-level improvements in multigrid for lattice QCD on large-scale computers}

\author[wuppertal]{Jesus Espinoza-Valverde}
\author[wuppertal]{Andreas Frommer}
\author[wuppertal]{Gustavo Ramirez-Hidalgo\corref{mycorrespondingauthor}}
\author[wuppertal,lausanne]{Matthias Rottmann}
\address[wuppertal]{Department of Mathematics, Bergische Universit\"at Wuppertal, 42097 Wuppertal, Germany}
\address[lausanne]{School of Computer and Communication Sciences, EPFL, 1015 Lausanne, Switzerland}
\cortext[mycorrespondingauthor]{Corresponding author (Email: g.ramirez@math.uni-wuppertal.de)}

\begin{abstract}
Numerical simulations of quantum chromodynamics (QCD) on a lattice require the frequent solution of linear systems of equations with large, sparse and typically ill-conditioned matrices. Algebraic multigrid methods are meanwhile the standard for these difficult solves. Although the linear systems at the coarsest level of the multigrid hierarchy are much smaller than the ones at the finest level, they can be severely ill-conditioned, thus affecting the scalability of the whole solver. In this paper, we investigate different novel ways to enhance the coarsest-level solver and demonstrate their potential using  DD-$\alpha$AMG, one of the publicly available algebraic multigrid solvers for lattice QCD. We do this for two lattice discretizations, namely clover-improved Wilson and twisted mass. For both the combination of two of the investigated enhancements, deflation and polynomial preconditioning, yield significant improvements in the regime of small mass parameters.
In the clover-improved Wilson case we observe a significantly improved insensitivity of the solver to conditioning, and for twisted mass we are able to get rid of a somewhat artificial increase of the twisted mass parameter on the coarsest level used so far to make the coarsest level solves converge more rapidly.
\end{abstract}


\end{frontmatter}

\nolinenumbers
\renewcommand\refname{}

\section{Introduction}
\label{sect:intro}

In lattice Quantum Chromodynamics (QCD) the Dirac operator describes the interaction between quarks and gluons in the framework of quantum field theory. Results of lattice QCD simulations represent essential input to several of the current and planned experiments in elementary particle physics (e.g., BELLE II, LHCb, EIC, PANDA, BES III \citep{alves2008lhcb,boca2014panda,belle2019belle,national2018assessment,yuan2019besiii}). Those simulations meanwhile achieve the required high precision (see e.g.\ \citep{borsanyi2021leading}) made possible through progress in both hardware and algorithms; see, e.g., \citep{joo2019status} and \citep{kahl2018adaptive}.

Inversions of discretizations of the Dirac operator are significant constituents in lattice QCD simulations. They arise in time critical computations such as when generating ensembles via the Hybrid Monte Carlo algorithm \citep{duane1987hybrid} or when computing observables \citep{Creutz_book,Gattringer_and_Lang_book} and they represent major supercomputer applications \citep{joo2019status}. 
There are multiple approaches to discretizing the continuum Dirac operator 
resulting in different discretized operators $D$, e.g.\ Wilson, Twisted Mass, Staggered, and others \citep{frezzotti2000local,Gattringer_and_Lang_book,kogut1975hamiltonian,wilson1974confinement}. 

Adaptive algebraic multigrid methods have established themselves as the most efficient methods for solving linear systems involving the Wilson-Dirac operator  \citep{babich_et_al_adapt_multigrid_for_latt_wilson,brannick_et_al_adapt_multigrid_for_lqcd,frommer_aggrg_multil_in_LQCD,adaptive_MG_in_LQCD_Rottmann,osborn_et_al_multigrid_for_clover}, and they have also found their way into other lattice QCD discretizations (see e.g.~\citep{alexandrou2016adaptive,brower2018multigrid,brower2020multigrid}). They demonstrate significant  speedups  compared to previously used conventional  Krylov  subspace solvers, achieving orders of magnitude faster convergence and showing insensitivity to the conditioning induced by the mass parameter. 

Multigrid methods use representations of the original operator on different levels $\ell = 1,\ldots,L$ with dimension decreasing with $\ell$.
They alternate steps of smoothing at each level with operations of transport between consecutive levels and (possibly quite inaccurate) solves of
systems with the coarsest operator on level $L$; see e.g.\ \citep{saad2003iterative}.

In adaptive algebraic multigrid implementations for lattice QCD, $L=2$ or $L=3$ levels are typically used \citep{luscher2012openqcd,clark2016accelerating} and the coarsest level is usually solved via a Krylov-based method, e.g.\ GMRES, possibly enhanced with a simple preconditioner and explicit deflation. In a parallel environment, the coarsest level solves suffer from an unfavorable ratio of communication vs.\ computation: A processor will have relatively few components to update, but a matrix vector multiplication will require a relatively high amount of data to be communicated between neighboring processors. Even more importantly, the communication cost for global reductions (mainly arising in the form of dot products) becomes very large as compared to computation \citep{Rottmann:2016gmm}.
As a result, coarsest level solves usually dominate the cost and routinely can reach up to 85\%  
of the overall computation time in, for example, current simulations with the twisted mass discretization \cite{Finkenrath2022b}. 
Hence, improving scalability of the coarsest-level 
solver is mandatory to improve the scalability of the whole multigrid solver. 

In this work, we consider a combination of three major techniques to improve the coarsest level solves:
\begin{enumerate}
\item Reducing the number of iterations by preconditioning with operators which do not require global communication.
\item Reducing the number of iterations by approximate deflation using Krylov recycling techniques.  
\item  Hiding global communication by rearranging computations.  
\end{enumerate}
As it will turn out, these techniques can yield a solver much less sensitive to conditioning when approaching the critical mass (see Section~\ref{sect:Wilson_results_latest_tuning_shifting_m0}) and improve scalability (see Section~\ref{sect:Wilson_results_latest_tuning_strong_scaling}). Furthermore, in the particular case of the twisted mass discretization, they prove to be useful in eliminating an artificially introduced parameter at the coarsest level (see Section~\ref{sect:twisted_mass_results_latest}).

The remainder of this paper is organized as follows. A brief introduction into lattice QCD discretizations and the DD-$\alpha$AMG library as one representative of algebraic multigrid methods for discretized Dirac operators 
is given in section \ref{sect:multig_in_LQCD}. Then, we introduce and discuss the three improvements for the coarsest level solves in section~\ref{sect:comm_reducing}.
Finally, in section~\ref{sect:num_tests} both algorithmic and computational tuning of the new parameters of the improved coarsest-level solver are performed, as well as strong scaling tests, followed by a re-tuning of an important parameter for twisted mass fermions solves.
All of the implementations for this work were done within the DD-$\alpha$AMG library for twisted mass fermions \citep{bacchio2019simulating}.

\section{Adaptive Aggregation based Algebraic Multigrid for LQCD}
\label{sect:multig_in_LQCD}

\subsection{Discretizations of the QCD Dirac operator}

The massive Dirac equation

\begin{equation}
\label{eq:Dirac_equation_cont}
\mathcal{D}\psi + m\psi = \eta,
\end{equation}
describes the dynamics of quarks through the interacting bosons of the strong force known as gluons. The quark fields $\psi$ and $\eta$ depend on each point $x$ of the 3+1 Euclidean space-time. The scalar parameter $m$ does not depend on $x$ and sets the mass of the quarks. The massless Dirac operator $\mathcal{D}$ 
\begin{equation}
\label{eq:Dirac_matrix_cont}
\mathcal{D} = \sum_{\mu=0}^{3} \gamma_{\mu} \otimes (\partial_{\mu} + A_{\mu}),
\end{equation}
where $\partial_{\mu} = {\partial}/{\partial x_{\mu}}$, encodes the gluonic interaction via the gluon (background) field  $A$ with $A_\mu$ from the Lie algebra $\mathfrak{su}(3)$. 
 The unitary and Hermitian matrices $\gamma_{\mu} \in \mathC^{4 \times 4}$ represent the generators of the Euclidean Clifford algebra 
\begin{equation}
\label{eq:Clifford_algebra}
\gamma_{\mu}\gamma_{\nu} + \gamma_{\nu}\gamma_{\mu} = 2\delta_{\mu\nu} I_{n}.
\end{equation}

Thus, at each point $x$ in space-time, the field $\psi(x)$ has 12 components corresponding to all possible combinations of 4 spins (acted upon by the $\gamma_\mu)$ and 3 colors (acted upon by the $A_\mu)$.  

The only known way to obtain predictions in QCD from first principles and non-pertubatively is to discretize space-time and then simulate on a computer. The discretization is typically formulated on a periodic $N_{t}{\times}N_{s}^{3}$ lattice $\mathcal{L}$ with uniform lattice spacing $a$; $N_{s}$ and $N_{t}$ represent the number of lattice points in the spatial and time directions, respectively.

In this paper we consider the two most commonly used discretizations, namely the (clover improved) \textit{Wilson} discretization and the \textit{twisted mass} discretization. 

The Wilson discretization is obtained by replacing the covariant derivatives by centralized covariant finite differences on the lattice. It contains an additional, second order finite difference stabilization term to correct for the doubling problem \citep{Gattringer_and_Lang_book}. 
The Wilson discretization yields a local operator $D$ in the sense that it represents a nearest-neighbor coupling on the lattice $\mathcal{L}$.

Introducing shift vectors $\hat{\mu} = (\hat{\mu}_{0},\hat{\mu}_{1},\hat{\mu}_{2},\hat{\mu}_{3})^{T} \in \mathbb{R}^{4}$ such that (with $\delta_{\mu\nu}$ the Kronecker delta)

\begin{equation}
\label{eq:mu_vectors}
\hat{\mu}_{\nu} = a \delta_{\mu\nu},
\end{equation}
i.e.\ vectors that direct us in each of the four space-time directions, the action of the Wilson-Dirac matrix $D$ on a discrete quark field $\psi(x), x \in \mathcal{L},$ can be written as

\begin{equation}
\label{eq:Dirac_matrix_discr_wilson}
\begin{split}
(D\psi)(x) = \left( \frac{m_0+4}{a} \right) \psi(x) - \frac{1}{2a} \sum_{\mu=0}^{3} ( (I_{4}-\gamma_{\mu}) \otimes U_{\mu}(x) ) \psi(x+\hat{\mu}) \\ - \frac{1}{2a} \sum_{\mu=0}^{3} ( (I_{4}+\gamma_{\mu}) \otimes U_{\mu}^{H}(x-\hat{\mu}) ) \psi(x-\hat{\mu}).
\end{split}
\end{equation}
The \textit{gauge link} matrices $U_{\mu}(x) = e^{-aA_{\mu}(x+\frac{1}{2}\hat{\mu})}$
live in the Lie group SU(3), since $A_{\mu}(x) \in \mathfrak{su}(3)$. The lattice indices $x \pm \hat{\mu}$ are to be understood periodically, and the mass $m_{0}$ sets the quark mass on the lattice \citep{qft_on_latt_montvay}.

The clover-improved Wilson-Dirac discretization replaces the first term as
\begin{eqnarray}
\lefteqn{\left(\frac{m_0+4}{a}\right)\psi(x)} & & \nonumber \\
& \to & \frac{1}{a} \left( (m_{0}+4)I_{12} - \frac{c_{sw}}{32} \sum_{\mu,\nu=0}^{3} (\gamma_{\mu}\gamma_{\nu}) \otimes (Q_{\mu\nu}(x) - Q_{\nu\mu}(x)) \right) \psi(x) ,
\label{eq:Dirac_matrix_discr_wilson_clover}
\end{eqnarray}
where
$ 
Q_{\mu\nu}(x) = Q_{x}^{\mu,\nu} + Q_{x}^{\mu,-\nu} + Q_{x}^{-\mu,\nu} + Q_{x}^{-\mu,-\nu}
$ 
are sums of the \textit{plaquettes}
\begin{equation}
\label{eq:expression_plaquettes}
Q_{x}^{\mu,\nu} = U_{\nu}(x)U_{\mu}(x+\hat{\nu})U_{\nu}^{H}(x+\hat{\mu})U_{\mu}^{H}(x).
\end{equation}

The Sheikholeslami-Wohlert parameter $c_{sw}$ is tuned such that it removes $\mathcal{O}(a)$ errors in the covariant finite difference discretization; see  \citep{lattice_methods_for_qcd,Gattringer_and_Lang_book,improved_cont_limit_latt_action}, e.g.

The (clover improved) twisted mass Dirac discretization arises by first performing a change of basis in the continuum relying on the invariance under chiral rotations and then discretizing on the lattice \citep{frezzotti2000local}. It is given as
\begin{equation}
\label{eq:twisted_mass_op_matrix}
D_{TM}(\mu) = D + i \mu \Gamma_{5},
\end{equation}
where $D$ is the clover-improved Wilson discretization \eqref{eq:Dirac_matrix_discr_wilson}, \eqref{eq:Dirac_matrix_discr_wilson_clover}, $\mu > 0$ is the \textit{twisted mass} parameter and $\Gamma_{5}$ describes the multiplication with $\gamma_{5} = i \gamma_{1} \gamma_{2} \gamma_{3} \gamma_{4}$ on the spins on each lattice site, i.e.\ $\Gamma_5 
= I_{n_{\mathcal{L}}} \otimes \gamma_{5} \otimes I_3 $ with $n_{\mathcal{L}}=N_{t}N_{s}^{3}$.

\subsection{Aggregation-based algebraic multigrid}
\label{sect:ddalphaamg}

This work deals exclusively with improvements of the coarsest level solves in adaptive algebraic multigrid methods.
In the description of algebraic multigrid to follow we will therefore dwell on details only when it comes to the coarsest level system. 

All multigrid methods for discretizations of the Dirac matrix are aggregation based, so the systems $D_\ell$ on the coarser levels represent again a nearest neighbor coupling on an (increasingly smaller) lattice. Aggregates correspond to blocks of contiguous lattice points on one level and are fused into one site of the lattice on the next. Prolongation and restriction are orthogonal operators obtained by breaking a number $n_{tv}$ of test vectors over the aggregates. 
The test vectors approximate small eigenmodes and are obtained within a bootstrap setup phase. Prolongation and restriction can be arranged to preserve a two-spin structure on all levels, as is done in \citep{babich_et_al_adapt_multigrid_for_latt_wilson,brannick_et_al_adapt_multigrid_for_lqcd,adaptive_MG_in_LQCD_Rottmann,osborn_et_al_multigrid_for_clover}, in which case the number of variables per grid point becomes $2n_{tv}$. If ones does not aim at preserving spin structure, as in \citep{luescher_local_coherence}, e.g., the number of variables is just $n_{tv}$. In principle, $n_{tv}$ can depend on the level $\ell$, but typically $n_{tv} \approx 20$ is chosen to be constant over the levels. 

The smoothers used on the different levels are typically either the Schwarz alternating procedure (SAP) as in DD-$\alpha$AMG or OpenQCD \citep{adaptive_MG_in_LQCD_Rottmann,luescher_SAP_original,luescher_LQCD_with_DD} or a plain Krylov subspace method such as GMRES as in \citep{babich_et_al_adapt_multigrid_for_latt_wilson,brannick_et_al_adapt_multigrid_for_lqcd, osborn_et_al_multigrid_for_clover}. An advantage of SAP is that it is rich in local computations requiring only comparably few  (nearest neighbor) communication.   

On the coarsest level, if we order odd lattice sites before even lattice sites, the matrix $D_L$ has the structure
\begin{equation} \label{eq:coarsest_level_structure}
D_L =\begin{pmatrix} D_{oo} & D_{oe}\\
   D_{eo} & D_{ee}\end{pmatrix}
\end{equation}
Herein, $D_{eo}$ and $D_{oe}$ represent the coupling with the nearest neighbors on the lattice. For the standard Wilson discretization, the diagonal matrices $D_{ee}$ and $D_{oo}$ are just multiples of the identity. When we include the clover term, $D_{ee}$ and $D_{oo}$ describe a self-coupling between the variables at each lattice point, i.e.\ they are block diagonal with the size of each block equal to $n_{tv}$, the number of variables per grid points. If we take the spin-preserving approach, the self-coupling is between variables of the same spin only, i.e.\ we actually have two diagonal blocks, each again with size $n_{tv}$ which now is half the number of variables per lattice site.

In order to solve coarsest level systems 
\[
D_L \psi =\eta \Longleftrightarrow 
\begin{pmatrix} D_{oo} & D_{oe}\\
   D_{eo} & D_{ee}\end{pmatrix} \begin{pmatrix} \psi_o \\ \psi_e \end{pmatrix} = \begin{pmatrix} \eta_o \\ \eta_e \end{pmatrix},
\]
the standard approach is to solve the odd-even reduced system 
\[
(D_{ee} - D_{eo}D_{oo}^{-1}D_{oe})\psi_e = \eta_e - D_{eo} D_{oo}^{-1}\eta_o
\]
for $\psi_e$ with a Krylov subspace method like GMRES, possibly enhanced with a deflation procedure, and then retrieve $\psi_o = D_{oo}^{-1}(\eta_o - D_{oe}\psi_e)$.

For future reference we denote
\begin{equation} \label{eq:factored_form}
D_c = D_{ee} - D_{eo}D_{oo}^{-1}D_{oe}
\end{equation}
the odd-even reduced matrix of the system at the coarsest level. For each lattice site, it describes a coupling with the variables from the 48 even lattice sites at distance 2. This is why $D_c$ is not formed explicitly -- a matrix-vector multiplication with $D_c$ is rather done by using the factored form \eqref{eq:factored_form}, where each multiplication with $D_{eo}$ and $D_{oe}$ involves the 8 nearest neighbor sites only. Since $D_{oo}$ and $D_{ee}$ are diagonal across the the lattice sites, multiplying with them does not involve other lattice sites, and the explicit computation of their inverses boils down to compute inverses of the local, small $n_{tv} \times n_{tv}$ diagonal blocks they are made up from.   

\section{Improving the coarsest level solves}
\label{sect:comm_reducing}

\subsection{GMRES}

The Generalized Minimal Residual Method GMRES \citep{saad1986gmres,saad2003iterative} is the best possible method to solve 
\[
A x = b,
\]
which relies on matrix-vector multiplications only, in the sense that it takes its iterates $x_m$ from $x_0 + K_m(A,r_0)$ with the 
Krylov subspace 
\begin{equation}
    \mathcal K_m(A,r_0) = \text{span} \{r_0, Ar_0, 
    \dots,A^{m-1} r_0 \}, \text{ where } r_0 = b-Ax_0,
    \label{eq:krylov_def}
\end{equation}
in such a manner that the residual $b-Ax_m$ has minimal 2-norm. The GMRES method relies on the Arnoldi process which computes an orthogonal basis $v_1,\ldots,v_m$ of  $\mathcal K_m(A,r_0) $ as described in Algorithm \ref{alg:arnoldi_process}

\begin{algorithm}[H]
    \KwData{(residual) vector $r_0$, matrix $A$, dimension $m$} 
    \KwResult{orthonormal matrix $V_{m} = [v_1| \cdots | v_m] \in \mathC^{n \times m}$, and Hessenberg matrix $\overline{H}_{m} = (h_{ij}) \in \mathC^{(m+1) \times m}$.}
$\beta = \|r_{0}\|_{2}$;\\
$v_{1} = r_{0}/\beta$; \\
   \For{$j = 1, \ldots, m$}{
        \For{$i = 1, \ldots, j$}{
            $h_{i,j} = (Av_j,v_{i})$; 
        }
        $w_j = Av_j - \sum_{i=1}^j h_{i,j}v_i$; \\
        $h_{j+1,j} = \| w_{j} \|_{2}$;  \\
        \If{$h_{j+1,j}=0$}
             {STOP} 
        $v_{j+1} = w_{j} / h_{j+1,j}$;
    }
    \caption{Arnoldi process}
    \label{alg:arnoldi_process}
\end{algorithm}

The GMRES iterate $x_m$ is then obtained as 
\[
x_m = x_0 + V_mz_m, 
\]
where  $z_m \in \mathC^m$ solves the least squares problem 
\begin{eqnarray*}
z_m &=& \mbox{argmin}_{z \in \mathC^m} \|b-A(x_0+V_m z )\|   = \mbox{argmin}_{z \in \mathC^m} \| \beta e_1 - \overline{H}_m z \| , \\  & & e_1 \in \mathC^{m+1} \mbox{ the first unit vector}.
\end{eqnarray*}

Lines 4-6 in the Arnoldi process orthogonalize $Av_j$ against $v_1,\ldots,v_j$. This is done with the classical Gram-Schmidt procedure. It is known that the following
mathematically equivalent modified Gram-Schmidt procedure 
which uses the partially orthogonalized vector $w_j$ for the computation of the $h_{i,j}$ is numerically more stable. 

\noindent
\LinesNotNumbered
\begin{algorithm}[H]
 \For{$j = 1, \ldots, m$}{
         $w_j = Av_j$; \\
        \For{$i = 1, \ldots, j$}{
            $h_{i,j} = (w_{j},v_{i})$; 
            $w_{j} = w_{j} - h_{i,j} v_{i}$;
        }
}
\end{algorithm}

When using classical Gram-Schmidt as in Algorithm~\ref{alg:arnoldi_process}, on a parallel computer the computation of the $k$ inner products can be fused into one single global reduction operation. This represents a substantial saving when the computational work load per process is small, as it is typically the case when solving systems with the coarsest grid marix $D_c$. On the other hand, stability is a minor concern, since typically very inaccurate 
solves on the coarsest level are sufficient, reducing the initial residual by just one order of magnitude.
Even then, however, the number of iterations required is high (some hundreds), and this makes GMRES increasingly inefficient since the orthogonalizations in the Arnoldi process require $j$ vector-vector operations in step $j$, eventually becoming far more costly than the matrix-vector multiplication. This is why, usually, GMRES must be restarted, i.e.\ a dimension $m$ for GMRES is fixed, and once it is reached, a new cycle of GMRES is started taking the last GMRES iterate as the new initial guess. In such restarted GMRES, termed GMRES($m$), the iterates now do not satisfy a global minimal residual condition any more, so the number of iterations to achieve a given accuracy must be expected to increase, but this is outweighed by a lower average cost per iteration.    

In the context of restarted GMRES, \textit{preconditioning} becomes useful not only to reduce the chances of slow-downs due to restarting, but also to improve scalability by trading global reductions against increased local computations. In section~\ref{sect:preconditioning} we describe how we use a left block diagonal and a right polynomial preconditioner to accomplish this.

Furthermore, we explore the interplay of those preconditioners with a deflation+recycling method, specifically GCRO-DR \citep{parks2006recycling}, in section~\ref{sec:deflation_via_GCRODR}. Both preconditioners allow for a reduction of the iteration count in coarsest-level solves with the immediate advantage of having less dot products at that level, but at the expense of an increase in local work (i.e.\ matrix-vector multiplications), particularly for the polynomial preconditioner.
Moreover, preconditioning tends to cluster the spectrum of the preconditioned matrix such that the number of small eigenmodes becomes small. It is particular in such situations that deflation 
achieved using GCRO-DR can yield substantial accelerations of convergence, since the deflation will effectively remove all these small eigenmodes.
A minor computational drawback of using deflation+recycling is that recycling vectors have to be constructed, which implies more extra work, and they are then deflated via projection in each Arnoldi iteration. But those deflations can be merged with the Arnoldi dot products when using classical Gram Schmidt, which increases the risk of losing stability. However, as explained before, the coarsest-level is solved quite inaccurately.

Finally, in Section~\ref{sect:comm_hiding} we integrate pipelining with the algorithms obtained so far as a means  to hide global communications and thus further improve scalability of the coarsest-level solves.

\subsection{Preconditioning}
\label{sect:preconditioning}

Left preconditioning uses a non-singular matrix $B$ to transform the original system $D_c x= b$ into
\begin{equation}
\label{eq:left_preconditioning_with_blockJac}
\begin{split}
B^{-1}D_c x = B^{-1}b.
\end{split}
\end{equation}

GMRES is now performed for this system. This means that in each iteration we now have a multiplication with $B^{-1}D_c$, which is done as two consecutive matrix-vector multiplications. This is why multiplicaton with $B^{-1}$ should be easy and cheap in computational cost. We took $B$ to be an approximate block Jacobi preconditioner. More precisely $B$ is the bock diagonal of $D_{ee}$, where 
each $n_{tv} \times n_{tv}$ diagonal block $B_i$ corresponds to all variables associated with lattice site $i$.
We compute $B_i^{-1}$ once for each $i$, and then perform the matrix vector products with $B_i^{-1}$ as direct matrix vector products. This is computationally more efficient than computing an $LU$-factorization of $B_i$ and then perform two triangular solves each time we need to multily with $B_i^{-1}$. Note that  all multiplications with $B_i^{-1}$ can be done in parallel and they do not require any communication if we---as we always do---keep all variables for a given lattice site on one processor.  

The diagonal blocks of $D_{ee}$ are not identical to those of $D_c$, and one might expect to obtain a more efficient preconditioner when using those of $D_c$.
Computing the part of $D_{eo} D_{oo}^{-1} D_{oe}$ contributing to each diagonal block of $D_c$ can, in principle, be done in parallel, but it requires inversions of $D_{oo}$ and communication with neighboring processes due to the couplings present in $D_{eo}$ and $D_{oe}$.
Limited numerical experiments suggest that the additional benefit of incorporating this part into the block diagonal preconditioner is moderate, so that we used the more simple-to-compute preconditioner which works exclusively with the block diagonal of $D_{ee}$. Further numerical experiments (see Section~\ref{sect:num_tests}) indicate that this block preconditioner typically gives a reduction by a factor of roughly 1.5 in the iteration count, 
at very little extra computational effort.
It has no effect if there is no clover term, since then the diagonal blocks are all multiples of the identity.

\subsection{Polynomial preconditioner}
\label{sect:polyprec}

Right polynomial preconditioning for the system $D_c x = b$ amounts to fixing a polynomial $q$ such that $q(D_c)$ is an approximation to $D_c^{-1}$ and then solving
\[
D_cq(D_c)y = b
\]
for $y$, using GMRES while keeping track of the back transformation $x_k = q(D_c)y_k$. This implies that we have
\[
x_k \in x_0 + K_{k}(D_cq(D_c),r_0),
\]
and $x_k$ is such that the residual $r_k = b-D_cx_k$ is minimized over the space $x_0 + K_{k}(D_cq(D_c),r_0)$. If $q$ has degree $d-1$, we invest $kd$ matrix-vector multiplications to build the orthonormal basis of $K_{k}(D_cq(D_c),r_0)$, a space of dimension $k$, in the Arnoldi process. With the same
effort in matrix-vector multiplications, we can as well build the $kd$-dimensional subspace $K_{kd}(D_c,r_0)$ which contains $K_{k}(D_cq(D_c),r_0)$. This shows that for the same amount of matrix-vector multiplications the residual of the $k$-th GMRES iterate of the polynomially preconditioned system can never be smaller than that of the $kd$-th iterate of standard GMRES. 
In other words, polynomial preconditioning reduces the iteration count while possibly increasing the total number of matrix-vector products.
Polynomial preconditioning can nevertheless be attractive for mainly two reasons: 
The first is that the above trade-off can be reversed in the presence of restarts. 
Standard GMRES is slowed down due to restarts, so if the polynomial preconditioner is good enough to allow to perform preconditioned GMRES without restarts, we might end up with less matrix-vector multiplications. The second is that cost beyond the matrix-vector products  can become  dominant. At iteration $k$, $k$ inner products and $k$ vector updates are performed in the Arnoldi process. In a parallel computing framework, these inner products, which require global communication, may become dominant for already quite small values of $k$, while very often matrix-vector products exhibit better potential for parallelization and display more promising scalability profiles compared to inner products. This even holds if we use the less stable standard Gram-Schmidt orthogonalization within Arnoldi which allows to fuse all $k$ inner products into one global reduction operation as explained after Algorithm~\ref{alg:arnoldi_process}.

We give an indicative example: Assume that $q$ has degree 3 and that we need 100 iterations with polynomially preconditioned GMRES. This amounts to 400 matrix-vector multiplications and 100 global reduction operations (for fused inner products). Assume further that with standard GMRES we need 200 iterations, i.e.\ 200 matrix vector multiplications and 200 global reductions. Then the additional matrix-vector multiplications in polynomial preconditioning are more than compensated for by savings in global reductions as soon as those take more time than two matrix-vector multiplications. Furthermore, the polynomially preconditioned method also saves on the vector operations related to the orthogonalization.     

Several types of polynomial preconditioners have been suggested in the literature, based on Neumann series, Chebyshev polynomials or least squares polynomials; see \citep{saad1985practical}, e.g. These approaches typically rely on detailed \textit{a priori} information on the spectrum of the matrix.
Recently, \citep{embree2021polynomial,loe2019new}, based on an idea from \citep{reichel1992}, showed that polynomial preconditioning with a polynomial obtained from a preliminary GMRES iteration can yield tremendous gains in efficiency when computing eigenpairs. We will use their way of adaptively constructing the polynomial for the preconditioner as we explain now. 

In standard GMRES, an iterate $x_d \in x_0 + K_d(D_c,r_0)$ can be expressed as $x_0 + q_{d-1}(D_c)r_0$ with $q_{d-1}$ a polynomial of degree $d-1$, and thus $r_d = b-D_cx_d = (I-D_cq_{d-1}(D_c))r_0 =:p_d(D_c)r_0$. Since $r_d$ is made as small as possible in GMRES, we can consider the polynomial $q_{d-1}$ to yield a good approximation $q_{d-1}(D_c)$ to $D_c^{-1}$. Strictly speaking, this interpretation only holds as far as the action on the vector $r_0$ is concerned, but we might expect this to hold for the action on just any vector if $r_0$ is not too special like, e.g., a vector with random components.

As is explained in \citep{loe2019new,morgan1991computing}, the polynomial $q_{d-1}$ can be constructed from 
the harmonic Ritz values of $D_c$ with respect to the Krylov subspace $K_{d}(D_c,r_0)$. These are the eigenvalues $\theta_i$ of the matrix 
\begin{equation}
    (H_d + h^2_{d+1,d} f e^T_d) ,
\end{equation}
with $H_d$ and $h_{d+1,d}$ from the Arnoldi process, Algorithm~\ref{alg:arnoldi_process}, and $f = H_d^{-H}e_d$, see \citep{morgan2006harmonic}. The  polynomial $p_d(t) = 1-tq_{d-1}(t)$ with $r_d = p_d(D_c)r_0$ is then given as
\begin{equation}
    p_d(t) = \prod^d_{i=1}\left(1 - \frac{t}{ \theta_i}
    \right),
    \label{eq:gmres_pol_loe}
\end{equation}
and since $q_{d-1}$ interpolates the values $\frac{1}{\theta_i}$ at the nodes $\theta_i$, \eqref{eq:gmres_pol_loe} gives, after some algebraic manipulation, a representation for $q_{d-1}$ similar to the Newton interpolation polynomial formula as
\begin{equation}
    q_{d-1}(t) 
    =   \sum^{d}_{i=1} \frac{1}{\theta_i} \prod_{j=1}^{i-1} \left(1 - \frac{t}{
    \theta_j}\right).
    \label{eq:q}
\end{equation}
Here, by convention, the empty product is $1$.

With this representation we can compute  $q_{k-1}(D_c)v$ using $d-1$ matrix-vector products by summing over accumulations of multiplications with $I-\frac{1}{\theta_j}D_c$.

\subsubsection*{Leja ordering}

The representation \eqref{eq:q} depends on the numbering which we choose for the harmonic Ritz values $\theta_i$ while, of course, mathematically $q_d$ is independent of the ordering. In numerical computation, however, the ordering matters due to different sensitivities to round-off errors. If $D_c$ is not very well conditioned, the range of $1/\theta_i$ may cover several orders of magnitudes so that it is important not to have all the big or all the small values appear in succession; see \citep{reichel1992}.     
An ordering choice that works well for a wide variety of cases is the Leja ordering \citep{reichel2003}. An algorithm to Leja order harmonic Ritz values is given as Algorithm~\ref{alg:leja}.

\vspace{0.5cm}

\begin{algorithm}[H]

    \KwData{Set $\mathbb K = \{  \theta_k\}_{k=1}^{d}$ of harmonic
    Ritz values.}

    \KwResult{Set $\{ \theta^L_k \}_{k=1}^{d}$ Leja ordered harmonic
    Ritz values.}

Choose $\theta^L_1 \in \mathbb{K}$ such that $|\theta^L_1| = \max \{ |\theta|: \theta \in \mathbb{K}\}$

    \For{$k = 2, \cdots, d$}{
        Remove $\theta^L_{k-1}$ from $\mathbb{K}$
        
        Determine $ \theta^L_k \in \mathbb K$, such that;
       \begin{equation}
             \prod_{j=1}^{k-1} | \theta^L_k -  \theta^L_j
             | = \max_{ \theta \in \mathbb K} \prod_{j=1}^{k-1} |
             \theta - \theta^L_j|. 
        \end{equation}
    }
    \caption{Leja ordering of harmonic Ritz values}
    \label{alg:leja}
\end{algorithm}
\vspace{0.5cm}
Algorithm~\ref{alg:pp_gmres} summarizes the process of computing and applying the preconditioning polynomial $q$ of degree $d-1$.

\vspace{0.5cm}

\begin{algorithm}[H]

    Construct the decomposition $D_cV_d = V_{d+1} \bar H_d$ 
    by running $d$ steps of the Arnoldi process using a random
    initial vector $v_0$;

    Compute the harmonic Ritz values $ \theta_k$ of $D_c$ as the eigenvalues of $H_d + h^2_{d+1,d} f e^T_d$;

    Leja order the obtained harmonic Ritz values $ \theta_k$;

    Run GMRES($m$) using the right preconditioner
    $q(D_c) =  \sum^{d}_{i=1} \frac{1}{\theta_i} \left(I - \frac{1}{
    \theta_1} D_c \right) \cdots \left(I - \frac{1}{
    \theta_{i-1}} D_c \right).$

    \caption{Polynomialy Preconditioned GMRES($m$) }
    \label{alg:pp_gmres}
\end{algorithm}

\vspace{0.5cm}

For the new implementations in DD-$\alpha$AMG involved in this work
we merge the block diagonal (left) preconditioner with polynomial preconditioning. This means that the preconditioning polynomial $q$ is constructed using $B^{-1}D_c$ rater than $D_c$ so that the 
overall preconditioned system takes the form
\begin{equation}
    \label{eq:system_with_both_precs}
    \begin{split}
        (B^{-1}D_{c}) q(B^{-1}D_{c}) y &= B^{-1}b, \\
            x &= q(B^{-1}D_{c})y.
    \end{split}
\end{equation}

\subsection{Deflation with GCRO-DR}
\label{sec:deflation_via_GCRODR}

A typical algebraic multigrid solve will take 10--30 iterations. Depending on the multigrid cycling strategy, each iteration will require one or more approximate solves on the coarsest level. In DD-$\alpha$AMG as in other multigrid methods for lattice QCD, the cycling strategy is to use K-cycles \citep{Notay2007}. This means that already with just three levels we will typically have in the order of ten coarsest-level solves per iteration, and this number increases as the number of levels increases. Even when using the preconditioning techniques presented so far, it usually happens that we need to perform several cycles of restarted GMRES to achieve the (relatively low) target accuracy required for these solves. 

This is why acceleration through deflation appears as an attractive approach in our situation. The idea is to use information gathered in one cycle of restarted GMRES to obtain increasingly better approximations of small eigenmodes of $D_c$ and to use those to augment the 
Krylov subspace for the next cycle. Then, even better approximations are computed and used in the following cycle or for the next system solve, etc. Effectively, this means that small eigenmodes are (approximately) deflated from the residuals, thus resulting in substantial acceleration of convergence.

Many deflation and augmentation techniques have been developed in the last 20 years \citep{coulaud2013deflation,soodhalter2020survey}, and some of them have already been used in lattice QCD, in particular eigCG \citep{stathopoulos2010computing} in the Hermitian case and GMRES-DR \citep{morgan2002gmres} for non-Hermitian problems. For Hermitian matrices, eigCG mimics the approximation of eigenpairs as done in Lanczos-based eigensolvers with a restart on the eigensolving part of the algorithm but in principle no restarts of CG itself. GMRES-DR, on the other hand, is used for non-Hermitian problems and, similarly to eigCG, it deflates approximations of low modes. When a sequence of linear systems is to be solved, eigCG can be employed in the Hermitian case \rev{and GMRES-DR/GMRES-Proj \citep{morgan2002gmres} can by applied in the general case. Here} we suggest to use GCRO-DR, the  generalized conjugate residual method with inner orthogonalization and deflated restarts of \citep{parks2006recycling}. GCRO-DR is partially based on GMRES-DR \citep{morgan2002gmres} and GCROT \citep{de1999truncation}, and it is particularly well suited to our situation where we not only have to perform restarts but also repeatedly solve linear systems with the same matrix and different right hand sides. 

A high level description of one cycle of GCRO-DR with a subspace dimension of $m$ is as follows\footnote{For the complete step-by-step GCRO-DR algorithm, see the appendix in \citep{parks2006recycling}.}:
\begin{enumerate}
    \item Extract $k<m$ approximations to small eigenmodes of $D_c$ from the current cycle.
    \item Determine a basis $u_1,\ldots.u_k$ for the space spanned by these approximate eigenmodes,
          gathered as colmns of $U_k = [u_1 | \cdots | u_k]$ such that $C_k = D_cU_k$ has orthonormal columns $c_1,\ldots,c_k$.
    \item With the current residual, perform $m-k$ steps of the Arnoldi process orthogonalizing not only against the newly computed Arnoldi vectors, but also against $c_1,\ldots,c_k$. This yields the relation (with $C_{k} = D_{c}U_{k}$)
    \begin{align} \label{eq:gcrodr_Arn_relation}
D_{c}[ U_{k} \ V_{m-k} ] &= [ C_{k} \ V_{m-k+1} ] G_m 
\end{align}
where  
\begin{align} G_m &= \begin{bmatrix}
I_{k} & B_{k}\\
0 & \bar{H}_{m-k}
\end{bmatrix} \text{ with $B_{k} = C_{k}^{H} D_{c} V_{m-k}$}
\end{align}
\item Obtain the new iterate by requiring the norm of its residual to be minimal over the space spanned by the columns of $U_k$ and $V_{m-k}$. This amounts to solving a least squares problem with $G_m  \in \mathC^{(m+1)\times m}$.
\end{enumerate}
The very first cycle of GCRO-DR, where we do not yet have approximate eigenmodes available, is just a standard GMRES cycle of length $m$. In all subsequent cycles, including those for solving systems with further right hand sides, the first step above typically computes the small eigenmodes as the small harmonic Ritz vectors of the matrix $G_m$ of the previous cycle. We have the option to stop updating the small approximate eigenmodes once they are sufficiently accurate. 

A source of extra work in GCRO-DR compared to GMRES is the construction of the recycling vectors, i.e.\ the columns in $U_k$ and $C_k$, from the harmonic Ritz vectors. By the imposition of $C_{k} = D_{c}U_{k}$ and via a QR decomposition of a very small matrix (which, in a parallel environment, is done redundantly in each process without the need for communications), $C_{k}$ can be efficiently updated from $U_{k}$, so this extra work is relatively small.
Another source of extra work is the deflation of $C_{k}$ in each Arnoldi iteration, which is again relatively cheap as the dot products due to those deflations can be merged with the already existing dot products from Arnoldi into one single global reduction.

\subsection{Communication hiding: pipelining}
\label{sect:comm_hiding}

An additional way to reduce communication time in the coarsest-level solves, complementary to what we have discussed so far, is \textit{communication hiding}, i.e.\ the overlapping of global communication phases with local computations. In this direction, \citep{ghysels2013hiding} introduces a \textit{pipelined} version of the GMRES algorithm which loosens the data dependency between the application of the matrix vector products and the dot products within the Arnoldi process. This is achieved by lagging the generation of the data obtained from the computation of the matrix-vector products from its actual use in the classical Gram-Schmidt process, so that the proposal vector that gets orthogonalized in a given iteration is
precomputed one or more iterations in the past. The number of ``lagging'' iterations is called the latency of the pipelining.

The method (with latency 1) requires the introduction of an extra set of vectors $v^a_{i}$ to store the matrix-vector products which are computed in advance.  These ``precomputed'' vectors are related to the orthonormal vectors $v_i$ from the Arnoldi process as $v^a_0 = v_0$ and $v^a_j = (A-\sigma_j I)v_j$ for $j \geq 1$. Here, $\sigma_j \in \mathC$ can be chosen arbitrarily, and judicious choices contribute to maintain numerical stability. In our particular implementations we have chosen the coefficients $\sigma_j$ to be zero since stability never was problematic due to the low relative tolerance required for the coarsest-level solves.

While this approach allows to hide global communications in Gram-Schmidt orthogonalizations behind the matrix-vector multiplication which typically requires only local communication, it does not \rev{do} so for 
the global communication required when normalizing the thus orthogonalized vector $w_j$ to obtain the final orthonormal vector $v_{j+1}$. 

As was observed in \citep{ghysels2013hiding}, this communication can be avoided if instead we compute $\|A v_j\|_2^2$ and then use
\begin{equation}
    h_{j+1,j}^2 = \|w_{j} \|^{2} = \| A v_{j} \|^{2}  - \sum_{i=1}^j |h_{i,j}|^2.
\end{equation}
The computation of $\|A v_j\|^2_2$ can now be done within the same global reduction communication used for the inner products yielding the coefficients $h_{i,j}$ in the classical Gram-Schmidt orthogonalization. 

It should be mentioned that a possible complication of this approach is that, due to numerical loss of orthogonality, we might get that 
$\| A v_{j} \|_{2}^{2}  - \sum_{i=1}^j |h_{i,j}|^2$ is not positive, which results in a breakdown. Even when there is no such breakdown, the above re-arrangement of terms tends to make the method  numerically less stable. 

\begin{algorithm}[h]
\SetAlgoLined
	$v_0 \leftarrow r_0 / \| r_0 \|_{2}$ \\	$v^a_0 \leftarrow B^{-1}D_c M v_0$, \mbox{where $M$ is the polynomial preconditioner $q_{d-1}(B^{-1}D_{c})$,} cf. sect. \ref{sect:polyprec}, with $B$ the block diagonal preconditioner, cf. sect. \ref{sect:preconditioning} \\
	\For{$i = 1, \ldots, m$}
	{
	$w_{i-1} \leftarrow B^{-1}D_c M v^a_{i-1}$ \\
	\For{$j=0,\ldots,i-1$} 
	{
	$h_{j,i-1}  \leftarrow (v_j, v^a_{i-1})$ 
	} 	
	$t \leftarrow \|v^a_{i-1}\|_{2}^2 - \sum_{j=0}^{i-1} h_{j,i-1}^2$ \\
	$\textbf{if} \hspace{1.2mm} t<0 \hspace{1mm} \textbf{then} \hspace{1mm} \text{breakdown}$ \\
	$h_{i,i-1} \leftarrow \sqrt{t}$ \\	
	$ v_i \leftarrow \left( v^a_{i-1} - \sum_{j=0}^{i-1} v_j h_{j,i-1} \right) / h_{i,i-1}$ \\
	$ v^a_i \leftarrow \left( w_{i-1} - \sum_{j=0}^{i-1} v^a_j h_{j,i-1} \right) / h_{i,i-1}$ 
	}
	$y \leftarrow \text{argmin}_{z \in \mathC^m} \big\| H_{m+1,m} z - \| r_0 \| e_1 \big\|$ \\
	$x\leftarrow x_0 + MV_m y$
	\caption{Latency 1 pipelined preconditioned GMRES}
    \label{alg:pipe1fgmres}
\end{algorithm}

Algorithm \ref{alg:pipe1fgmres} summarizes the method.
Mathematically, it is equivalent to standard GMRES. As compared to the latter, pipelined GMRES
requires more memory to store the extra set of vectors $v_i^a$, and it requires more local computation in the form of additional AXPY operations.
The advantage is that the global communications required to obtain the $h_{j,i-1}$ coefficients in the outer loop $i$ and $\|v_{i-1}^a\|^2$
can be performed in parallel with the matrix-vector multiplication needed to compute $w_{i-1}$.

When extending the pipelining approach to GCRO-DR and thus \rev{including} recycling and deflation, a new set of vectors to be computed in advance must be introduced to store the matrix-vector plus preconditioner application on the recycling vectors.
\rev{The resulting pipelined GCRO-DR algorithm seems not to have been proposed in the literature so far. We formulate it as Algorithm~\ref{alg:pipe1gcrodr}, where we refer to specific lines in the original description of GCRO-DR in the Appendix of \cite{parks2006recycling} for those parts which are identical but somewhat lengthy to reproduce. Note that the columns of the matrix $V^{c} = B^{-1}D_{c}M C_{k}$ need to be updated whenever $U_{k}$ (and by association $C_{k}$) changes, which, in Algorithm~\ref{alg:pipe1gcrodr} is done on every iteration of the while loop.
Although updating $U_{k}$ and $C_{k}$ remains cheap, as in the original GCRO-DR algorithm, constructing $V^{c}$ is now expensive. This is why we will introduce later in section \ref{sect:Wilson_results_latest_tuning} a modification where we stop updating the recycling subspace once this has been a given number of times.}

\begin{algorithm}[h]
	\SetAlgoLined
        \rev{initial steps as done in lines 1-19 of the original GCRO-DR algorithm. This produces $U_{k}$, $C_{k}, x$ and $r_{1}$ \\
        $r \leftarrow r_{1}$} \\
        \While{\rev{$\|r\|_{2} > tol$}}{
	$v_0 \leftarrow r / \| r \|$ \\
        $\textbf{for} \hspace{1mm} j=0,\ldots,k-1  \hspace{1.2mm} \textbf{do} \hspace{1.5mm}  v^c_{j}= \rev{B^{-1}}D_{c} M c_j$,  \hspace{2mm} 
        \\
        \hspace*{1.5ex} where $M$ is the polynomial preconditioner $q_{d-1}(B^{-1}D_{c})$ (cf.\ \\
        \hspace*{1.5ex} sect.\ \ref{sect:polyprec}) with $B$ the block diagonal preconditioner from sect. \ref{sect:preconditioning} \\	   $v^a_0 \leftarrow \rev{B^{-1}} D_{c} M v_0$ \\
	\For{$i = 1, \ldots, m-k$}
	{
		$w \leftarrow \rev{B^{-1}} D_{c} M v^a_{i-1}$ \\
		$\textbf{for} \hspace{1mm} j=0,\cdots,\rev{k-1}  \hspace{1.2mm} \textbf{do} \hspace{1.5mm} b_{j,i-1} = (c_j,v^a_{i-1})$   
         \\
		$\textbf{for} \hspace{1mm} j=0,\cdots,i-1  \hspace{1.2mm} \textbf{do} \hspace{1.5mm} h_{j,i-1} = (v_j,v^a_{i-1})$  
        \\
		$t \leftarrow \|v^a_{i-1}\|_{2}^2 - \sum_{j=0}^{i-1} h_{j,i-1}^2 - \sum_{j=0}^{i-1} b_{j,i-1}^2$ \\
		$\textbf{if} \hspace{1.2mm} t<0 \hspace{1mm} \textbf{then} \hspace{1mm} \text{breakdown}$ \\
		$h_{i,i-1} \leftarrow \sqrt{t}$ \\	
		$ v_i \leftarrow \left( v^a_{i-1} - \sum_{j=0}^{i-1} v_j h_{j,i-1} 
		- \sum_{j=0}^{i-1} c_j b_{j,i-1} \right) / h_{i,i-1}$ \\
		$ v^a_i \leftarrow \left( w - \sum_{j=0}^{i-1} v^a_j h_{j,i-1} - \sum_{j=0}^{i-1} v^c_j b_{j,i-1} \right) / h_{i,i-1}$ \\
	}
        \rev{define $D_{k}$, $\widetilde{U}_{k}$, $\widehat{V}_{m}$, $\widehat{W}_{m+1}$ and $G_{m}$ as in lines 23-26 of the original  GCRO-DR. What we denote $G_{m}$ corresponds to $\underbar{G}_{m}$ in \citep{parks2006recycling}} \\
	$y \leftarrow \text{argmin}_{z \in \mathC^m} \big\| G_{m}z - \|r\|_{2} e_{k+1} \big\|$ \\
	\rev{$x\leftarrow x + M \widehat{V}_{m}y$, $r \leftarrow r - \widehat{W}_{m+1} G_{m} y$} \\
        \rev{update $U_{k}, C_{k}$ as in lines 30-34  of the original GCRO-DR algorithm}
        }
	\caption{Latency 1 pipelined preconditioned GCRO-DR}
	\label{alg:pipe1gcrodr}
\end{algorithm}

\subsection{The whole coarsest-level solver}

We have combined the block diagonal preconditioning with $D_{ee}^{-1}$, adaptive polynomial preconditioning, deflation and recycling via GCRO-DR, and pipelining into a single implementation for coarsest-level solves within DD-$\alpha$AMG. The code was implemented in a modularized way such that the user can enable any combination of these four methods already during the compilation of the program.
The implementation is ready to use and available at a GitHub repository\footnote{ \href{https://github.com/JesusEV/DDalphaAMG_ci}{https://github.com/JesusEV/DDalphaAMG\_ci}}. At runtime, \rev{the coarsest level is adaptive in the following sense: the setup phase of DD-$\alpha$AMG, as described in \citep{frommer2013adaptive}, changes the matrices at coarser levels as it improves the multigrid hierarchy. Hence whenever such a change happens for the matrix at the coarsest level, a set of flags are set such that the next time that the coarsest-level solver is called, the block-diagonal preconditioner is reconstructed, then the polynomial is reconstructed, and finally the recycling subspace of GCRO-DR is reconstructed. The construction of the polynomial is done by calling Arnoldi with a random coarsest-level right hand side as the first Arnoldi vector. In the case of GCRO-DR, the inner GMRES (step 3 in the high level description of GCRO-DR in section~\ref{sec:deflation_via_GCRODR}) is forced to undergo an Arnoldi of length (at least) $k$ when the recycling data $U_{k}$ and $C_{k}$ need to be reconstructed}.

Extensive numerical tests of all the methods implemented are presented in Section~\ref{sect:num_tests}.

\section{Numerical tests}
\label{sect:num_tests}

All computations were performed on the JUWELS supercomputer from the J\"ulich Supercomputing Centre. In most of our tests on JUWELS, one process per node and 48 OpenMP threads\footnote{The number of threads per MPI process can be varied in some cases and for example a value of 20 might be better in certain situations where the number of nodes is extremely large and the work per thread is very small such that the thread barriers start becoming significant. We take this into consideration in section~\ref{sect:Wilson_results_latest_pipelining}.} per process were used in the JUWELS cluster module. In Section \ref{sect:Wilson_results_latest} we present results for the  clover-improved Wilson discretization using configuration D450r010n1 from the D450 ensemble of the CLS collaboration (see e.g.\ \citep{mohler2018cls}). Section \ref{sect:twisted_mass_results_latest} deals with twisted mass fermions where we use configuration conf.1000 of the cB211.072.64 ensemble of the Extended Twisted Mass Collaboration (see \citep{alexandrou2018simulating}). In both cases the lattice is of size $128 \times 64^3$.

\subsection{The clover-improved Wilson operator}
\label{sect:Wilson_results_latest}

With the inclusion of new algorithms at the coarsest level, new parameters for these algorithms need to be tuned, which we do in Section \ref{sect:Wilson_results_latest_tuning}.
Once this is done, we test how the whole solver is affected by the numerical conditioning of the discretized operator by varying the mass parameter. Moreover, we perform some scaling tests of the whole solver.

The block diagonal preconditioner of section \ref{sect:preconditioning} is always used in all experiments here. It comes at very little extra computational cost, but can give a reduction by a factor of up to about 1.5 in the iteration count at the coarsest level, as can be seen from Table~\ref{tab:effect_of_BDP}. 

\begin{table}
\caption{Effect of the block diagonal preconditioner (BDP) on coarsest-level solves in DD-$\alpha$AMG, where we have the BDP as the only preconditioner of GMRES. The second and third columns are average number of iterations at the coarsest level in the solve phase. We have used configuration D450r010n1 here with different values of $m_{0}$.
}
\centering
 \begin{tabular}{||c c c||} 
 \hline
 $m_{0}$ & without BDP & with BDP \\
 \hline\hline
 -0.3515 & 21 & 15 \\ 
 -0.35371847789 & 35 & 26 \\ 
 -0.354 & 40 & 28 \\ 
 -0.3545 & 52 & 37 \\ 
 \hline
 \end{tabular}
 \label{tab:effect_of_BDP}
\end{table}

\subsubsection{Tuning parameters}
\label{sect:Wilson_results_latest_tuning}

The set of default parameters in our solves can be found in Table~\ref{tab:ddalphaamg_default_params} (see \citep{frommer2013adaptive,Rottmann:2016gmm} for more on these parameters).

\begin{table}[h]
\caption{Base parameters in our DD-$\alpha$AMG solves.}
\centering
 \begin{tabular}{||c c c||} 
 \hline
 $\ell=1$ & restart length of FGMRES & 10 \\
 & relative residual tolerance & $10^{-9}$ \\
 & number of test vectors & 24 \\
 & size of lattice-blocks for aggregates & $4^4$ \\
 & pre-smoothing steps & 0 \\
 & post-smoothing steps & 3 \\
 & Minimal Residual iterations in SAP& 4 \\
 \hline
 $\ell=2$ & restart length of FGMRES & 5 \\
 & maximal restarts of FGMRES & 2 \\
 & relative residual tolerance & $10^{-1}$ \\
 & number of test vectors & 32 \\
 & size of lattice-blocks for aggregates & $2^4$ \\
 & pre-smoothing steps & 0 \\
 & post-smoothing steps & 2 \\
 & Minimal Residual iterations in SAP& 4 \\
 \hline
 $\ell=3$ & restart length of GMRES & 60 \\
 & maximal restarts of GMRES & 20 \\
 & relative residual tolerance & $10^{-1}$ \\
 \hline
 \end{tabular}
 \label{tab:ddalphaamg_default_params}
\end{table}

Next to this set of parameters, there are now two more that need to be tuned in order to minimize the total execution time of the solves. These parameters are:

\begin{itemize}
\item $k$: number of recycling vectors 
in GCRO-DR.
\item $d$: degree of the polynomial employed as polynomial preconditioner.
\item $u$: the number of times we update the recycling subspace information represented by $U$ and $C$. After $u$ updates, we continue to use the last $U$ and $C$ in all further restarts and all further solves with new right hand sides. 
\end{itemize}

The performance dependence on $u$ shows initial gains for smaller values of $u$ with only marginal progress for already only moderately large values. This is why we fixed $u=10$ in all our experiments.

\begin{figure}
\includegraphics[width=.5\textwidth]{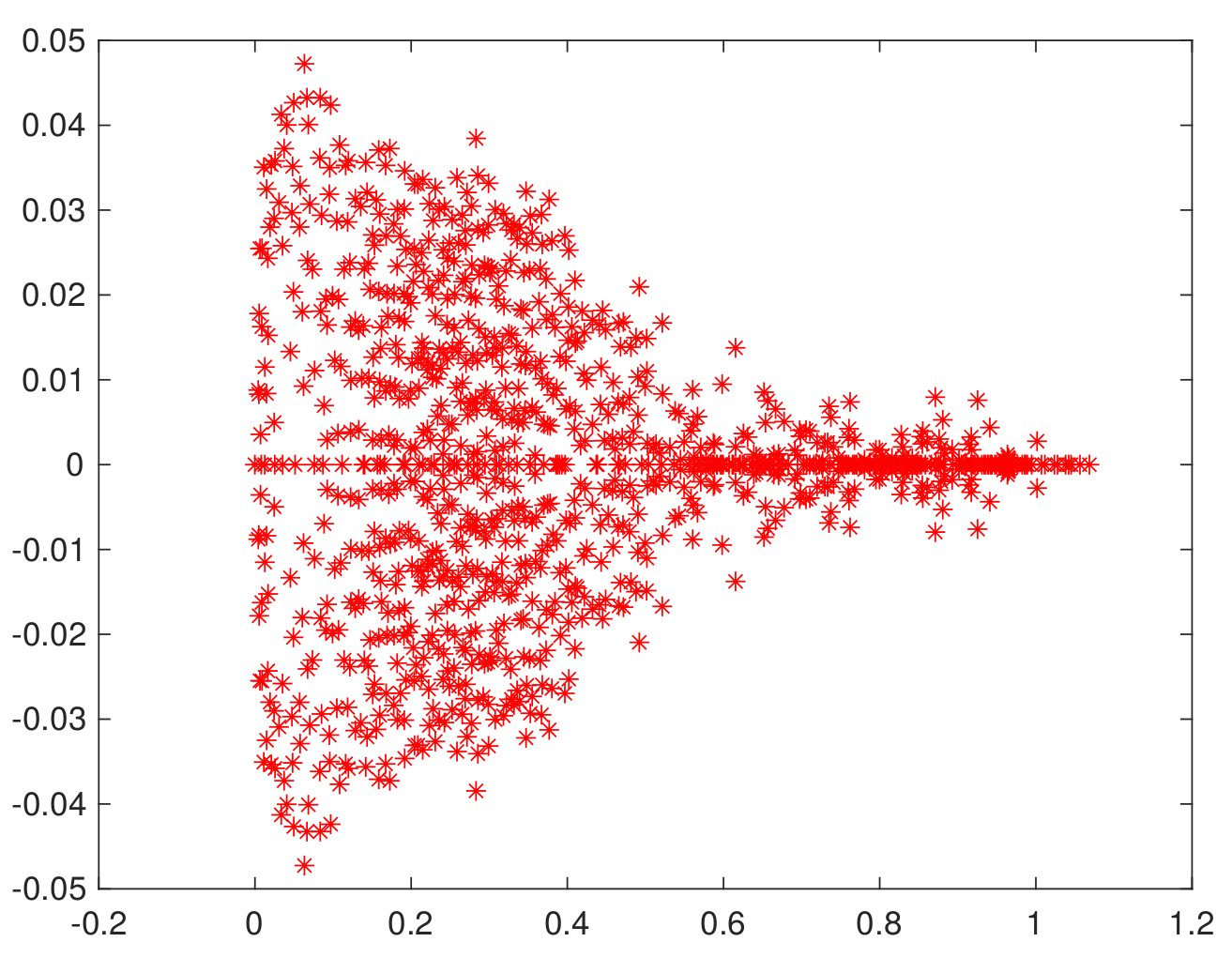} \hfill
\includegraphics[width=.48\textwidth]{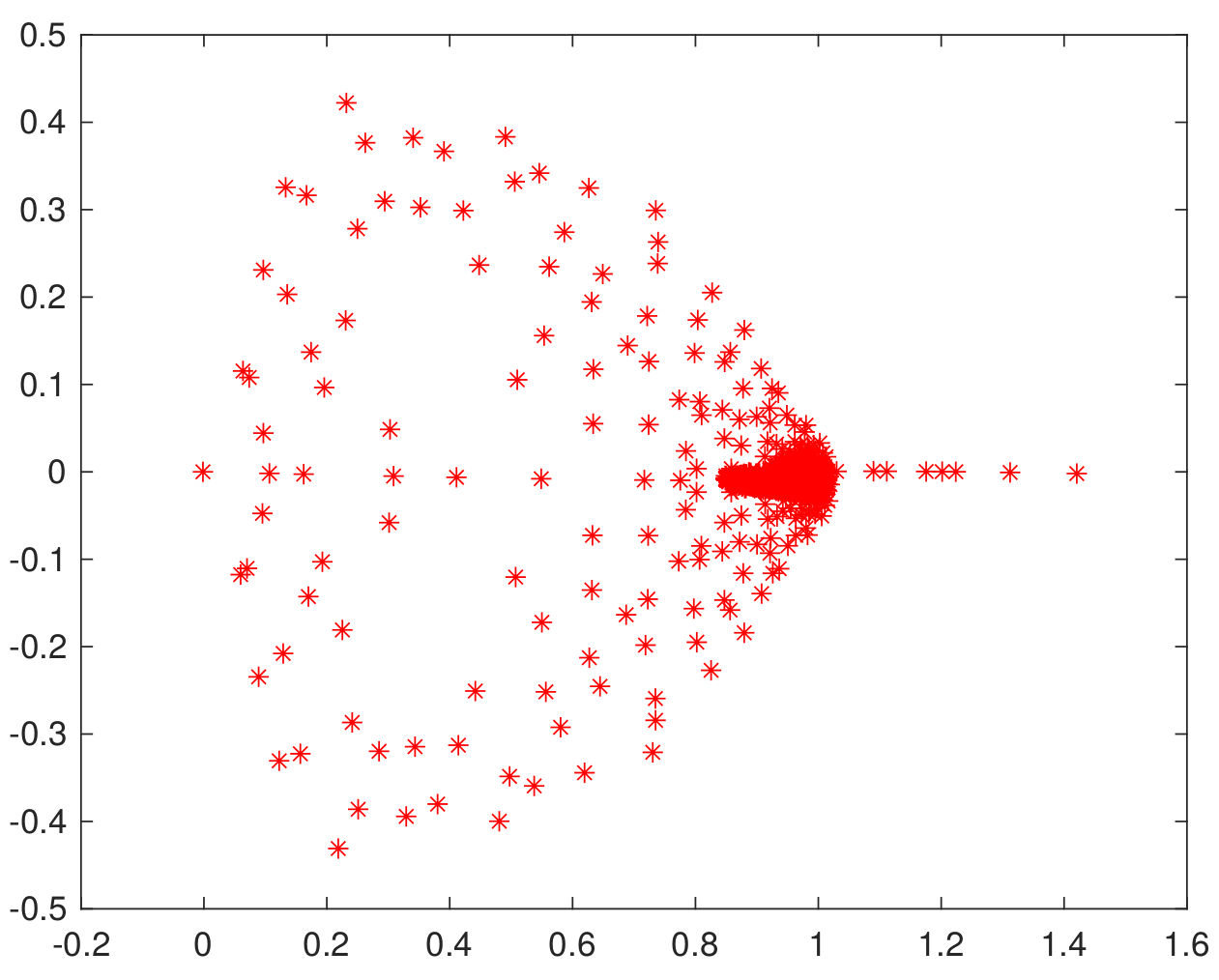}
\caption{Spectrum of $D_{c}$ (left) and of $D_{c} q(D_{c})$ (right) in a four-level hierarchy for a $64 \times 32^3$ clover-improved Wilson configuration. Due to odd-even preconditioning and the number of test vectors from levels three to four being 32, 
$D_{c}$ is of size $1024 \times 1024$. \label{fig:Wilson_estimates_exectime}}
\end{figure}

There is an interplay between the polynomial preconditioner and GCRO-DR, that needs to be kept in mind when tuning the new parameters. Preconditioning has a tendency to move most of the small eigenmodes close to 1, which is why we can remove the remaining small eigenmodes with relatively moderate values of $k$ in GCRO-DR. So preconditioning actually helps deflation to become more efficient. To illustrate this, we have constructed a four-level method for a clover-improved Wilson configuration\footnote{This configuration was provided by the lattice QCD group at the University of Regensburg via the Collaborative Research Centre SFB-TRR55, with parameters $m_{0} = -0.332159624413$ and $c_{sw} = 1.9192$ \citep{Bali:2014gha}. We shifted the diagonal of the Dirac matrix to $m_{0} = -0.33555$ to obtain a more ill-conditioned system.} of size $64 \times 32^3$. The aggregation at different levels is such that the coarsest-level lattice is of size $4 \times 2^3$. The spectrum of both $D_{c}$ and $D_{c} q(D_{c})$ are presented in Figure~\ref{fig:Wilson_estimates_exectime}, where we see that deflation (via e.g.\ GCRO-DR) would be of benefit due to the scattered nature of the low modes of $D_{c} q(D_{c})$.

\begin{figure}[h]
  \begin{minipage}[b]{0.4\textwidth}
    \hspace*{-1cm}
    \includegraphics[width=1.45\textwidth]{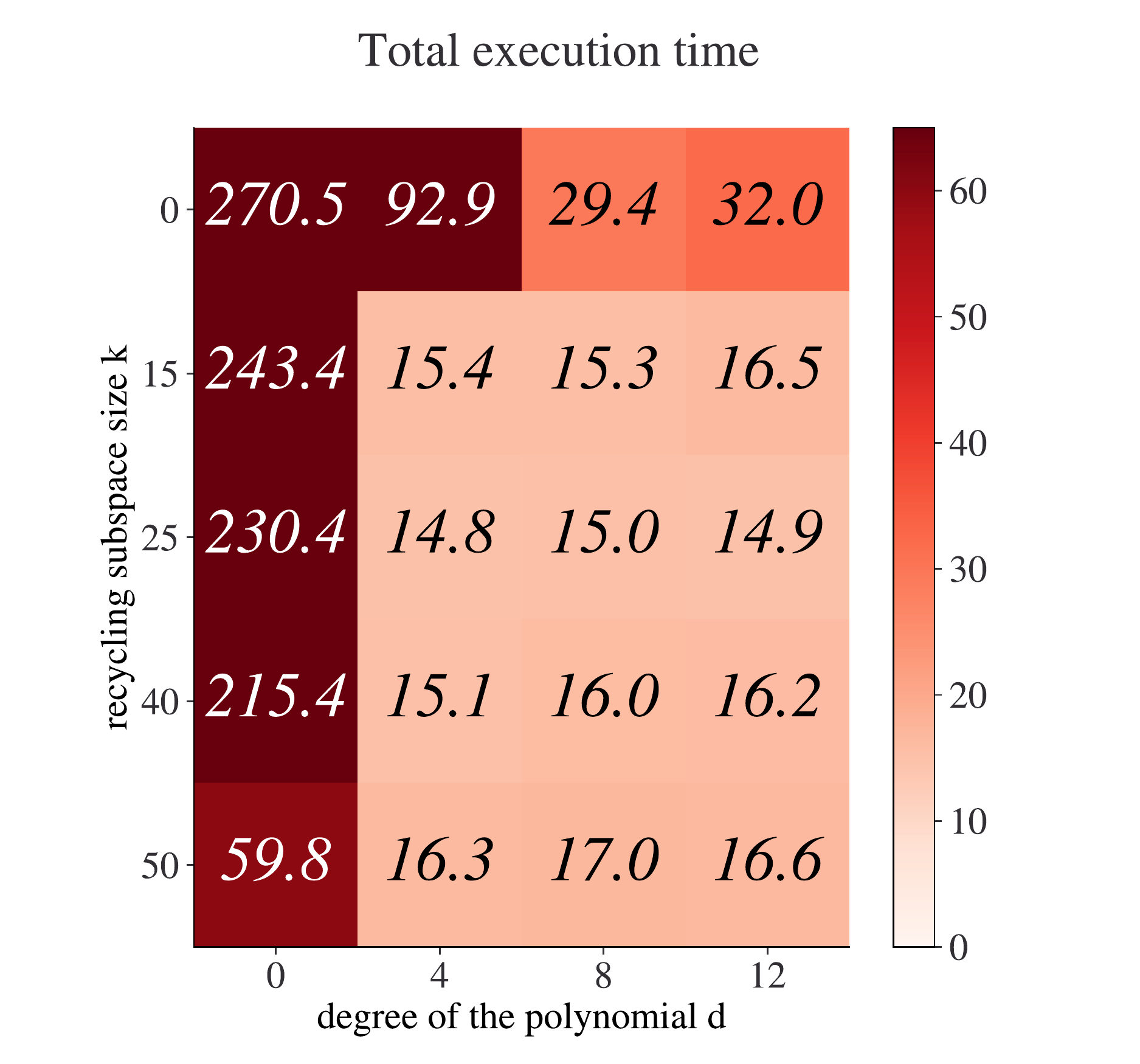}
    \hspace{2mm}
  \end{minipage}
  \hfill
  \begin{minipage}[b]{0.4\textwidth}
   \hspace*{-1cm}
   \includegraphics[width=1.45\textwidth]{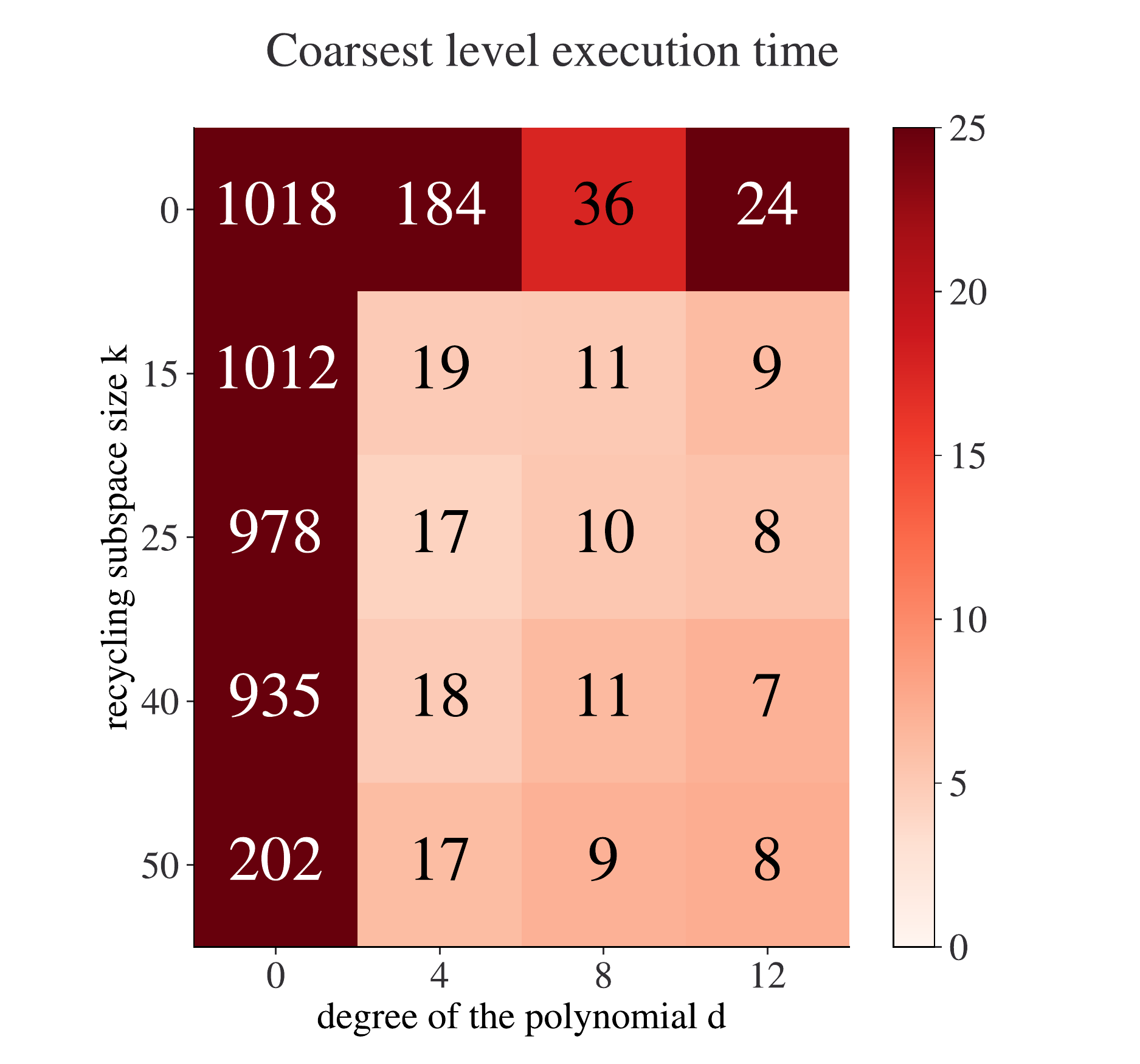}
   \hspace{2mm}
  \end{minipage}
  \caption{Tuning of the parameters $k$ and $d$ . The color of each square in the heatmap from the left represents the total execution time of the whole DD-$\alpha$AMG solver, while the right corresponds to the time spent at the coarsest level. The configuration was for a lattice of size $128{\times}64^{3}$; we used 32 nodes with 48 OpenMP threads, each. All these computations were done for $m_{0}=-0.355937$ (i.e. the most ill-conditioned case in Figure~\ref{fig:versusm0_Wilson}). The numbers in the boxes on the right indicate the average number of iterations at the coarsest level during the whole multigrid solve phase\rev{, and the numbers in the boxes on the left are the execution times of the whole solve phase}.}
  \label{fig:tuning_kd_Wilson}
\end{figure}

With a $128 \times 64^3$ lattice and with $u=10$ fixed, 
Figure~\ref{fig:tuning_kd_Wilson} displays the execution time for the solve phase in DD-$\alpha$AMG for a cartesian product of $(k,d)$ pairs.
The choice $k=0, d=0$ (left upper) corresponds to no polynomial \rev{preconditioning} and no deflation. We found that the choice $d=4, k=25$ gives more than a factor of 18 improvement over this case, and that choices for $k,d$ in the neighborhood of this optimal pair affect the execution time only marginally.
As a rule, smaller values of  $k$ should be preferred as a low value of $k$ reduces the risk of inducing instabilities (due to having deflation and Arnoldi dot products merged, see Section~\ref{sec:deflation_via_GCRODR}).
From Figure~\ref{fig:tuning_kd_Wilson} we see that $d=4$ and $d=8$ are equally good in the particular tests that we have run.

\subsubsection{Pipelining}
\label{sect:Wilson_results_latest_pipelining}

\begin{table}
\caption{Effect of pipelining on the whole DD$\alpha$AMG solver. We have used configuration D450r010n1 here with $m_{0} = -0.355937$. \rev{For GCRO-DR and the polynomial preconditioner, the pair $(k,d)=(25,4)$ was used.}
}
\centering
 \begin{tabular}{||c c c c||} 
 \hline
 $N_{proc}$ & $N_{thr}$ & with pipelining & without pipelining \\
 \hline\hline
 128 & 20 & 5.02 & 4.56 \\ 
 256 & 20 & 3.18 & 2.98 \\ 
 512 & 10 & 2.7 & 2.6 \\ 
 1024 & 10 & 2.18 & 2.0 \\ 
 \hline
 \end{tabular}
 \label{tab:pipelining_results_high_level}
\end{table}

Table~\ref{tab:pipelining_results_high_level} gives a comparison of the total execution time for one DD$\alpha$AMG-solve without and with pipelining in the preconditioned GCRO-DR solves on the coarsest level. As we see, pipelining always {\em increases} the execution time by up to 10\%, even for larger number of processors. \rev{For the numerical experiments of this section we have fixed $(k,d) = (25,4)$. Other $(k,d)$ pairs that we tested lead to similar results with the same  explanation holding in terms of the interplay of pipelining with the polynomial preconditioner for the specific surface-to-volume ratio that we deal with here at the coarsest level.}

\begin{table}
\caption{Execution times for parts of the \rev{coarsest} grid solves with and without pipelining. Times in last three columns are in seconds. \rev{For GCRO-DR and the polynomial preconditioner, the pair $(k,d)=(25,4)$ was used.} }
\centering
 \begin{tabular}{||c c c c c c||} 
 \hline
 $N_{proc}$ & $N_{thr}$ & pipel. & mvm & mvm-w. & glob-reds\\ 
 \hline\hline
 256 & 20 & OFF & 0.638 & 0.127 & 0.235 \\ 
 512 & 10 & OFF & 0.687 & 0.210 & 0.259 \\ 
 1024 & 10 & OFF & 0.606 & 0.144 & 0.33 \\ 
 256 & 20 & ON & 0.964 & 0.331 & 0.0558 \\ 
 512 & 10 & ON & 0.901 & 0.294 & 0.0614 \\ 
 1024 & 10 & ON & 1.050 & 0.468 & 0.0896 \\ 
 \hline
 \end{tabular}
 \label{tab:pipelining_results_low_level}
\end{table}

To understand this behavior better, we timed the relevant parts of the computation and communication on the coarsest level for our MPI implementation on JUWELS. The results are reported in Table~\ref{tab:pipelining_results_low_level}. Here, for different numbers $N_\mathit{proc}$ of processors and $N_\mathit{th}$ of threads, we report three different timings: mvm refers to the time spent in one matrix-vector multiplication (arithmetic plus communications), mvm-w is the time
processors spent in an MPI-wait for the communication related to the matrix-vector multiplication to be completed. 
Note that DD-$\alpha$AMG uses a technique from \citep{krieg2010tuning} that aims to overlap computation and communication as much as possible for the nearest neighbor communication arising in the matrix-vector multiplication. 
Finally, glob-reds reports the time processors wait for the global reductions to be completed. We see that these wait times are indeed almost entirely suppressed in the pipelined version. 
However, we also see that we do not succeed to hide the communication for the global reductions behind the matrix-vector multiplication, since the time of the latter is increased when pipelining is turned on.  
We conclude that on JUWELS and with the MPI implementation in use, the communication for global reductions and for the matrix-vector multiplication compete for the same network resources, thus counteracting the intended hiding of communication.         
\rev{One might hope that this clash of pipelining with the halo exchanges can be alleviated if} the matrix-vector multiplications present in the polynomial preconditioner can be done in a non-synchronized manner, so that the mvm-wait times are almost reduced to zero. \rev{This} can be achieved in a manner similar to what is called {\em communication avoiding} GMRES \citep{hoemmen2010communication}, \rev{specifically via the \textit{matrix powers kernels} \cite{powerkernels}}, whereby one exchanges vector components which belong to lattice sites up to a distance $k$ in one go and then can evaluate polynomials up to degree $k$ in $A$ without any further communication. \rev{It is important to note that, as parallelism is increased and the number of lattice sites per process is correspondingly reduced, distance $k$ halos
comprise many more lattice points than $k$ times those of a distance 1 halo. This means that the communication \textit{volume} is substantially increased, and this might become prohibitive.} \rev{Communication avoiding methods at the coarsest level have been used in the context of lattice QCD already, see e.g.\ \citep{ayyar2022optimizing}.}
Implementing \rev{them} is a major endeavor, though, and out of scope for this paper.

\subsubsection{Scaling with the conditioning}
\label{sect:Wilson_results_latest_tuning_shifting_m0}

\begin{figure}
\centering
\hspace*{-0.5cm} 
\includegraphics[width=1.1\textwidth]{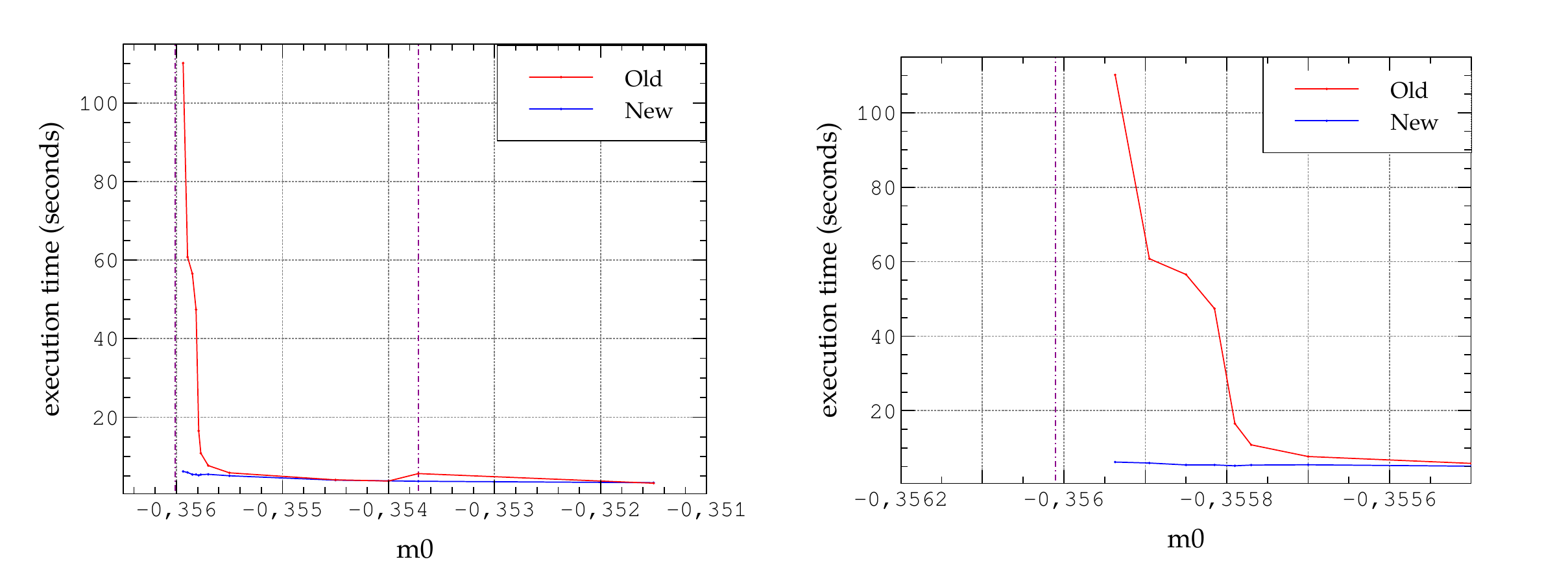}
\caption{Total execution time of the solve phase in DD-$\alpha$AMG as the system becomes more ill-conditioned (i.e.\ as $m_{0}$ becomes more negative). The vertical dashed line closest to -0.354 represents the value with which the Markov chain was generated and the vertical dashed line closest to -0.356 represents $m_\mathit{crit}$. The right plot zooms into the region where the old version of the solver does not perform well.}
\label{fig:versusm0_Wilson}
\end{figure}

With the number of processes fixed at 128, we now shift the mass parameter $m_{0}$ to see how much the new coarse grid solver improves upon the old when conditioning of the Wilson-Dirac matrix changes. Results are given in Figure~\ref{fig:versusm0_Wilson}. Pipelining is turned off for these tests and we put $k=25$, $d=4$ throughout. \rev{We decrease the mass parameter below the point which was used when generating the configuration in the HMC process such as to obtain quite ill-conditioned matrices. We stopped our experiments at mass values relatively close to $m_\mathit{crit}$, the first value for which the Wilson-Dirac matrix becomes singular. In the methods using deflation, the eigenmode corresponding to the smallest eigenvalue is effectively suppressed, which is why execution times do not vary significantly even if we moved further close to $m_\mathit{crit}$. This is in contrast to the non-deflated method, where the increasing ill-conditioning requires more and more iterations on the coarsest level.}

Our experiments show that the improvements due to the new coarse grid solver are close to marginal for the better conditioned matrices, but that they become very substantial in the most ill-conditioned cases.
Actually, the times for the new solver are almost constant over the whole range for $m_0$, whereas the times for the old one increase drastically for the most ill-conditioned systems. The left dotted vertical line, located very close to -0.356, represents the location of $m_{crit}$, i.e.\ the value of $m_{0}$ for which the Dirac operator becomes singular.
In Table~\ref{tab:mg_iters_with_m0} we summarily report an interesting observation regarding the setup of DD$\alpha$AMG. According to the bootstrap principle, the setup performs iterations in which the 
multigrid hierarchy is improved from one step to the next. In ill-conditioned situations, the solver on the coarsest level might stop at the prescribed maximum of possible iterations rather than because it has achieved the required accuracy. This affects the quality of the resulting final operator hierarchy. The table shows that for a given comparable effort for the setup, the one that uses the improved coarse grid solvers obtains a better overall method, since the coarsest systems are solved more accurately.
\begin{table}
\caption{Number of iterations of the outermost FGMRES in DD-$\alpha$AMG as $m_{0}$ moves down to more ill-conditioned cases.}
\centering
 \begin{tabular}{||c c c||} 
 \hline
 $m_{0}$ & old & new\\ [0.5ex] 
 \hline\hline
 $-0.355770$ & 19 & 19 \\ 
 $-0.355815$ & 23 & 19 \\
 $-0.355850$ & 26 & 20 \\
 $-0.355895$ & 29 & 20 \\
 $-0.355937$ & 30 & 20 \\ [1ex] 
 \hline
 \end{tabular}
 \label{tab:mg_iters_with_m0}
\end{table}

\subsubsection{Strong scaling}
\label{sect:Wilson_results_latest_tuning_strong_scaling}

\begin{figure}
\centering
\hspace*{-0.7cm}
\includegraphics[width=0.65\textwidth]{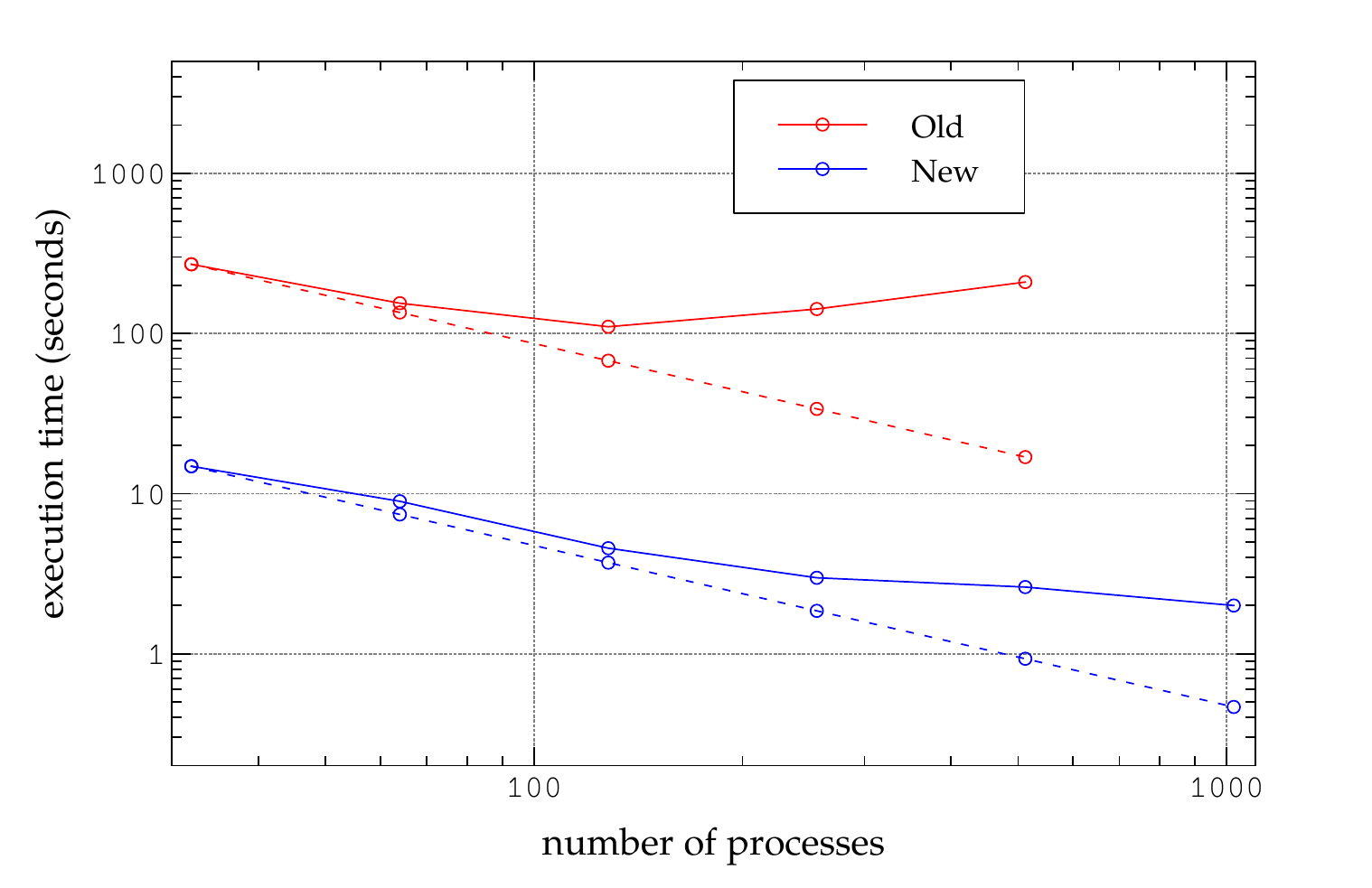}
\caption{Strong scaling tests on Wilson fermions for the new coarsest-level additions. The solves were applied over a $128{\times}64^{3}$ lattice. \textit{Old} means the previous version of DD-$\alpha$AMG without the coarsest-level improvements discussed in this paper, and the vertical axis represents the whole solve time.
The dashed lines indicate how both cases would behave in case of perfect scaling. All these computations were done for $m_{0}=-0.355937$ (i.e. the most ill-conditioned case in Figure~\ref{fig:versusm0_Wilson}).}
\label{fig:strong_scaling_Wilson}
\end{figure}

Figure~\ref{fig:strong_scaling_Wilson} reports a strong scaling study for both the old and the new coarse grid solves within DD-$\alpha$AMG. We see that the new version improves scalability quite substantially, due to the better scalability of the coarse grid solve. This is to be attributed to the fact that the new coarse grid solver reduces 
the fraction of work spent in inner products and thus global reductions which start to dominate for large numbers of processors.  Still, we do not see perfect scaling, one reason being that after some point the wait times occurring in the matrix-vector multiplications become perceptible.

We have mentioned before the need to tune the number of OpenMP threads when going to a very large number of processes. Indeed, as the work per OpenMP thread becomes quite small, then thread barriers start becoming problematic from a computational performance point of view. The results presented in Figure~\ref{fig:strong_scaling_Wilson} already include this tuning: for 32 and 64 process we used 48 OpenMP threads per process, for 128 and 256 we switched to 20 threads, and finally for 512 and 1024 we rather used 10 threads. Pipelining was kept off for these scaling tests.

\begin{figure}
\centering
\hspace*{-0.7cm}
\includegraphics[width=0.65\textwidth]{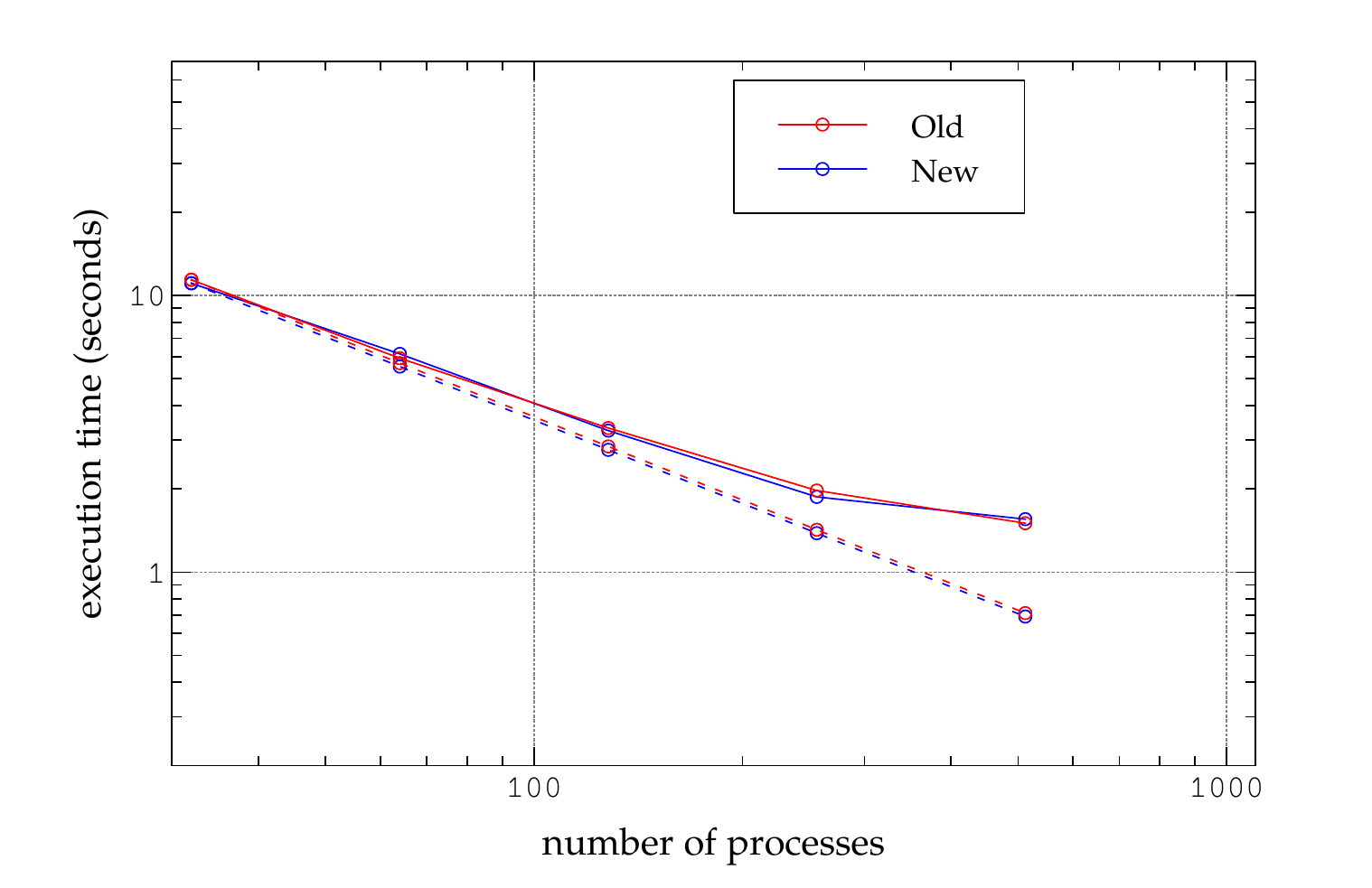}
\caption{Strong scaling tests on Wilson fermions for the new coarsest-level additions. The solves were applied over a $128{\times}64^{3}$ lattice. \textit{Old} means the previous version of DD-$\alpha$AMG without the coarsest-level improvements discussed in this paper, and the vertical axis represents the whole solve time.
The dashed lines indicate how both cases would behave in case of perfect scaling. All these computations were done for $m_{0}=-0.35371847789$.}
\label{fig:strong_scaling_Wilson_originalm0}
\end{figure}

We ran another strong scaling test with $m_{0} = -0.35371847789$, i.e.\ the value of the mass parameter originally used for the generation of the ensemble. The results of this are shown in Figure~\ref{fig:strong_scaling_Wilson_originalm0}. The coarsest-level improvements did not bring visible gains to the whole solver execution time, which is due to having a quite well-conditioned coarsest level. We also note that the new and old versions of the solver match in this case and for any other relatively well-conditioned value $m_{0}$, which gives consistency to the new implementation with respect to the old one.

\subsection{The twisted mass operator}
\label{sect:twisted_mass_results_latest}

We now turn to the twisted mass discretization
\eqref{eq:twisted_mass_op_matrix}, where the parameter $\mu$ ``shields" the spectrum away from 0 in the sense that the smallest singular value of $D_{TM}$ is $\sqrt{\lambda_{sm}^{2} + \mu^{2}}$ with $\lambda_{sm}$ the smallest eigenvalue in absolute value of the symmetrized clover-improved Wilson-Dirac operator $Q = \Gamma_{5} D$.
This is, in general, algorithmically advantageous, but eigenvalues now have the tendency to cluster around the smallest ones  \citep{frezzotti2000local}.

There is an extension of DD-$\alpha$AMG that operates on twisted mass fermions \citep{bacchio2019simulating}. In that version, the twisted mass parameter $\mu$ remains, in principle, propagated without changes from one level to the next. On the coarsest level, the clustering phenomenon of small eigenvalues is particularly pronounced, resulting in large iteration numbers of the solver at the coarsest level. A way to alleviate this is to use, instead of $\mu$, a multiple $\mu_{c} = \delta \cdot \mu$ with a factor $\delta>1$ on the coarsest level. As was shown in \citep{alexandrou2016adaptive}, this can decrease the required number of iterations substantially. 

We work with configuration conf.1000 of the cB211.072.64 ensemble of the Extended Twisted Mass Collaboration, see \citep{alexandrou2018simulating}. The lattice size is $128 \times 64^3$ and $\mu = 0.00072$.
Different values of $\mu_{c}$ lead to a different spectrum at the coarsest level, and therefore for each different value of $\mu_{c}$ a new tuning of the new coarsest-level parameters $u,k$ and $d$ has to be performed. We tuned parameters in a similar way as we described before for Wilson fermions.

\rev{For $\delta=8.0$, restarted GMRES with no further enhancements was used. Preconditioning and deflation were used with $\delta=1.0$, where we found that $d=20$ for the degree of the polynomial preconditioner and $k=400$ for the size of the recycling subspace in GCRO-DR lead to minimal execution times, with $u=10$ during the setup phase and $u=1$ during the solve phase. Pipelining was kept off during these tests due to the conclusions from Section \ref{sect:Wilson_results_latest_pipelining}. We ran strong scaling tests for these two values of $\mu_{c}$. Due to the different nature of the twisted mass discretization compared to the Wilson case, some parameters in the multigrid solver had to be changed with respect to Table \ref{tab:ddalphaamg_default_params}; we report these changes in Table \ref{tab:ddalphaamg_twisted_mass_params}. The value of 400 for the cycle length of restarted GMRES seems quite large, but this is needed only to ensure that enough memory is allocated to run an Arnoldi of length $k=400$ in every first construction of the recycling subspace by GCRO-DR. Furthermore, due to the coarsest level representing a significant part of the overall execution time, even when both preconditioners and GCRO-DR were under use, a careful tuning lead us to opt for 1 thread per process and 32 MPI processes per compute node.}

\begin{table}[h]
\caption{Parameters changed in DD-$\alpha$AMG, for twisted mass solves, with respect to the base ones in Table \ref{tab:ddalphaamg_default_params}.}
\centering
 \begin{tabular}{||c c c||} 
 \hline
 $\ell=1$ & relative residual tolerance & $10^{-10}$ \\
 & number of test vectors & 32 \\
 & post-smoothing steps & 4 \\
 \hline
 $\ell=2$ & post-smoothing steps & 3 \\
 \hline
 $\ell=3$ & restart length of GMRES & 400 \\
 & maximal restarts of GMRES & 10 \\
 \hline
 \end{tabular}
 \label{tab:ddalphaamg_twisted_mass_params}
\end{table}

\begin{figure}[h]
\centering
\hspace*{-0.8cm} 
\includegraphics[width=0.6\textwidth]{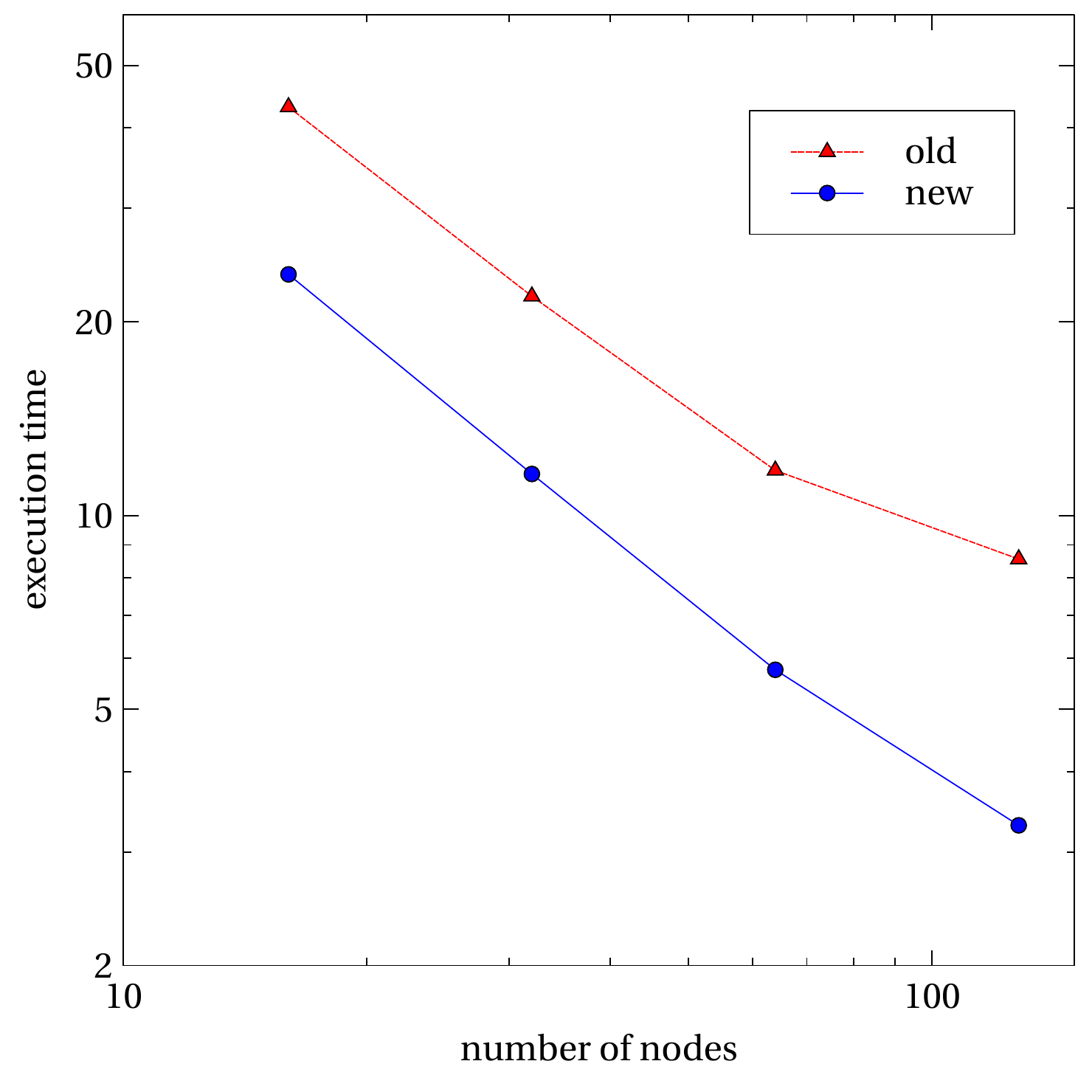}
\caption{Strong scaling tests on twisted mass fermions for the new coarsest-level additions. The total execution time of the solve phase is shown here.
\textit{Old}: $\delta = 8.0$, with the restarted GMRES as the coarsest-level solver with no further improvements. 
\textit{New}: $\delta = 1.0$, both preconditioners and GCRO-DR on, and with $d=20$ for the polynomial preconditioner and $k=400$ for the size of the recycling subspace in GCOR-DR. The number of nodes are 16, 32, 64 and 128, with 32 MPI processes per node and 1 thread per process. The smallest execution time for the old version is 8.6 seconds, and the corresponding one for the new version is 3.3 seconds.}
\label{fig:strong_scaling}
\end{figure}

\rev{The results of our strong scaling tests are shown in Figure~\ref{fig:strong_scaling} which shows the impact of the coarsest-level improvements on the overall performance of the solver. For 128 nodes, a speedup of 2.6 in the total execution time of the solve phase is obtained, and the scalability of the whole solver has improved. To achieve these gains, the polynomial preconditioner has been particularly crucial. The impact of the polynomial preconditioner on the coareset level solves 
is threefold. First, without polynomial precoditioning, 
the plain GMRES polynomial has difficulties in adjusting to the coarsest-level twisted mass spectrum when $\delta=1.0$. This means that we need 
very large (i.e.\ $>$ 2,000) GMRES restart lengths, with the iteration count to reach $10^{-1}$ still being of the order of 10,000 for many systems. Second, as illustrated in Figure \ref{fig:Wilson_estimates_exectime} for the Wilson case, the polynomial scatters the spectrum closest to zero, making it easier for GCRO-DR to resolve those small modes that need to be deflated. 
And lastly, the nearest-neighbor nature of the communications in the matrix-vector multiplications
in the polynomial preconditioner 
scale better than global communications, and furthermore polynomial preconditioning reduces the latter due to faster convergence.
}


\rev{Next to a considerable reduction in the total execution time, our approach brings $\delta$ back down to 1.0, i.e.\ $\mu_c = \mu$, thus removing this {\em ad hoc} and somewhat artificial parameter. Shifting $\delta$ down to 1.0 reduces the iteration counts of FGMRES at levels $\ell=1$ and $\ell=2$ from 20 and 7 (on average) to 16 and 3, respectively. This is because  we have
a more faithful representation of the correct value of $D_{c}$ when $\delta=1.0$, leading to better coarse-grid corrections at levels $\ell=2$ and $\ell=1$. A further benefit shows 
when taking
a closer look at the timings of the most demanding segments of the solver once the coarsest-level improvements have been included. We report these times in Table \ref{tab:timings_fastest_solver_state}. The times for the three components of the solver add up to 2.453s, which represents 74.3\% of the total execution time. Each of these times 
could be further reduced by using more advanced vectorization (we currently use SSE) and by storing the data at the coarsest level and for the smoothers in half precision instead of the currently used single precision.}

\begin{table}
\caption{Some execution times (in seconds) in the new version of the multigrid solver, applied to twisted mass, running on 128 nodes. The label mvm has the same meaning here as in Table \ref{tab:pipelining_results_low_level}.}
\centering
 \begin{tabular}{||c c||} 
 \hline
 function & exec.\ time\\ [0.5ex] 
 \hline\hline
 smoother $\ell=0$ & 0.874 \\ 
 smoother $\ell=1$ & 0.569 \\
 mvm $\ell=2$ & 1.01 \\ [1ex] 
 \hline
 \end{tabular}
 \label{tab:timings_fastest_solver_state}
\end{table}

\section{Conclusions and outlook}

We have extended the DD-$\alpha$AMG solver at its coarsest level, from a plain restarted GMRES to encompass four different methods: block-diagonal and polynomial preconditioning, GCRO-DR and pipelining. In the case of the clover-improved Wilson discretization, and for the particular matrix used here, the method recovers its (approximate) insensitivity to conditioning as we approach the critical mass. \rev{In the twisted mass case, we have reduced the overall execution time of the solve phase by a factor of approximately 2.6 and improved the scalability of the solver. At the same time,  we returned the twisted mass parameter $\mu_c$ on the coarsest level back to the more canonical value $\mu$ used on the other levels. The current state of the solver, when applied to twisted mass, suggests the need for better vectorization and the use of lower precision, which we will integrate in future versions of our code.}

An important outcome of the discussion in Section~\ref{sect:Wilson_results_latest_pipelining} is the need of a communication-avoiding scheme in our coarsest-level implementations: as we increase the number of nodes in our executions, nearest-neighbor communications become a two-fold problem, first in the sense of their lack of scalability, and second as they interfere with the global communications trying to be hidden by pipelining. We will implement a communication-avoiding method, which we expect to have a nice interplay with both pipelining and the polynomial preconditioner.

\section{Acknowledgments}

We would like to thank Francesco Knechtli and Tomasz Korzec for providing us with the clover-improved Wilson configuration files, and to Jacob Finkenrath and Simone Bacchio for providing the twisted mass ones. In particular, we would also like to thank Jacob for stimulating discussions regarding some simulations aspects and challenges particular to twisted mass. The authors gratefully acknowledge the Gauss Centre for Supercomputing e.V. (www.gauss-centre.eu) for funding this project by providing computing time through the John von Neumann Institute for Computing (NIC) on the GCS Supercomputer JUWELS at Jülich Supercomputing Centre (JSC), under the project with id CHWU29. This work was funded by the European Commission within the STIMULATE project.

\section{References}

\bibliography{references}

\end{document}